\documentclass{article} 
\usepackage[OT1]{fontenc}
\usepackage{etoolbox}
\newcommand{\arxiv}[1]{\iftoggle{colt}{}{#1}}
\newcommand{\colt}[1]{\iftoggle{colt}{#1}{}}
\newtoggle{colt}
\global\togglefalse{colt}

\arxiv{
\bibliographystyle{alpha}
\newcommand{\citep}[1]{\cite{#1}}
}

\colt{
\usepackage{times}
}
\usepackage[letterpaper, left=1in, right=1in, top=1in, bottom=1in]{geometry}

\usepackage{bbm}
\usepackage{bm}
\usepackage{multirow}
\usepackage{enumitem}
\usepackage{marginnote}
\usepackage{graphicx}
\usepackage{mdframed}
\usepackage{mathtools}
\usepackage{hyperref}       
\usepackage{microtype}      

\usepackage{amsfonts}       

\usepackage{colortbl}

\usepackage{todonotes}
\usepackage{amsthm}
\newtheorem{theorem}{Theorem}
\newtheorem{lemma}[theorem]{Lemma}

\newtheorem{definition}[theorem]{Definition}
\newtheorem{remark}[theorem]{Remark}

\newtheorem{assumption}{Assumption}
\newmdtheoremenv{condition}{Condition}
\AfterEndEnvironment{theorem}{\noindent\ignorespaces}
\AfterEndEnvironment{lemma}{\noindent\ignorespaces}
\AfterEndEnvironment{corollary}{\noindent\ignorespaces}
\AfterEndEnvironment{proposition}{\noindent\ignorespaces}
\AfterEndEnvironment{definition}{\noindent\ignorespaces}
\AfterEndEnvironment{assumption}{\noindent\ignorespaces}
\AfterEndEnvironment{proof}{\noindent\ignorespaces}
\usepackage{xcolor}

\renewcommand{\v}{\bm{v}}

\newcommand{\cB}{\mathcal{B}}

\newcommand{\wtilde}{\widetilde}

\newcommand{\loss}{\ell}
\DeclareBoldMathCommand{\vloss}{\loss}

\DeclareBoldMathCommand{\fakegrad}{\mathring{\bm{g}}}
\DeclareBoldMathCommand{\e}{e}
\DeclareBoldMathCommand{\p}{p}
\DeclareBoldMathCommand{\u}{u}
\DeclareBoldMathCommand{\w}{w}
\DeclareBoldMathCommand{\x}{x}
\DeclareBoldMathCommand{\l}{l}
\DeclareBoldMathCommand{\vzero}{0}
\let\top\intercal

\newcommand{\reals}{\mathbb{R}}

\DeclareMathOperator{\E}{\mathbb{E}}

\renewcommand{\x}{\bm{x}}   

\newcommand{\y}{\bm{y}} 
\newcommand{\g}{\bm{g}}

\DeclareMathOperator*{\argmin}{argmin}

\newcommand{\G}{\mathbf{G}}
\newcommand{\cK}{\mathcal{K}}
\newcommand{\bnabla}{\bm{\nabla}}
\newcommand{\inte}{\operatorname{int}}
\newcommand{\dom}{\operatorname{dom}}

\newcommand*{\what}[1]{\widehat{#1}}

\newcommand{\cA}{\mathcal{A}}

\newcommand{\nn}{\nonumber}

\newcommand{\inner}[2]{\langle#1,#2\rangle}

\newcommand{\bzeta}{\bm{\zeta}}
\newcommand{\veps}{\varepsilon}

\DeclareMathOperator{\sep}{{\normalfont{\texttt{sep}}}}
\DeclareMathOperator{\proj}{{\normalfont\texttt{proj}}}
\DeclareMathOperator{\grad}{{\normalfont{\texttt{grad}}}}

\newcommand{\Bfrak}{\mathfrak{B}}  
\newcommand{\bnu}{\bm{\nu}}  
\newcommand{\K}{\mathcal{C}}

\newcommand{\cE}{\mathcal{E}}

\newcommand{\vertiii}[1]{{\left\vert\kern-0.25ex\left\vert\kern-0.25ex\left\vert #1 
		\right\vert\kern-0.25ex\right\vert\kern-0.25ex\right\vert}}

\makeatletter
\newcommand*\bcdot{\mathpalette\bigcdot@{.5}}
\newcommand*\bigcdot@[2]{\mathbin{\vcenter{\hbox{\scalebox{#2}{$\m@th#1\bullet$}}}}}

\newlength\myindent
\usepackage{MnSymbol}
\usepackage{comment}
\usepackage{bbding}
\usepackage{booktabs}
\usepackage{footnote}
\usepackage[nameinlink,capitalize]{cleveref}

\let\norm\undefined

\newcommand{\norm}[1]{\left\lVert#1\right\rVert}

\makesavenoteenv{tabular}
\makesavenoteenv{table}
\usepackage{algorithm}
\usepackage{algorithmic}
\colt{
	\newcommand{\algcomment}[1]{  \textcolor{blue!70!black}{\scriptsize{\texttt{{//#1}}}}}
}
\arxiv{
	\newcommand{\algcomment}[1]{  \textcolor{blue!70!black}{\footnotesize{\texttt{{//#1}}}}}
}

\makeatletter
\g@addto@macro\bfseries{\boldmath}
\makeatother

\arxiv{
	\setlength{\marginparwidth}{10ex} 
	\setlength{\marginparsep}{5mm}
}

\colt{
\title[Efficient Online Exp-Concave Optimization]{Quasi-Newton Steps for Efficient Online Exp-Concave Optimization}

\coltauthor{%
	\Name{Zakaria Mhammedi} \Email{mhammedi@mit.edu}\\
	\AND
	\Name{Khashayar Gatmiry} \Email{gatmiry@mit.edu}\\
	\addr Massachusetts Institute of Technology%
}
}

\arxiv{
	\title{Quasi-Newton Steps for Efficient Online Exp-Concave Optimization}
\author{{\bf Zakaria Mhammedi}    \\
	{\bf Khashayar Gatmiry}  \\ Massachusetts Institute of Technology \\ \texttt{\{mhammedi, gatmiry\}@mit.edu}}
\date{November 23, 2022;\quad Revised: February 14, 2023}
}

\begin{document}

	\maketitle
	
	\begin{abstract}
		The aim of this paper is to design computationally-efficient and optimal algorithms for the online and stochastic exp-concave optimization settings. Typical algorithms for these settings, such as the Online Newton Step (ONS), can guarantee a $O(d\ln T)$ bound on their regret after $T$ rounds, where $d$ is the dimension of the feasible set. However, such algorithms perform so-called \emph{generalized projections} whenever their iterates step outside the feasible set. Such generalized projections require $\Omega(d^3)$ arithmetic operations even for simple sets such a Euclidean ball, making the total runtime of ONS of order $d^3 T$ after $T$ rounds, in the worst-case. In this paper, we side-step generalized projections by using a self-concordant barrier as a regularizer to compute the Newton steps. This ensures that the iterates are always within the feasible set without requiring projections. This approach still requires the computation of the inverse of the Hessian of the barrier at every step. However, using stability properties of the Newton iterates, we show that the inverse of the Hessians can be efficiently approximated via Taylor expansions for most rounds, resulting in a $\wtilde O(d^2 T +d^\omega \sqrt{T})$ total computational complexity, where $\omega\in(2,3]$ is the exponent of matrix multiplication. In the stochastic setting, we show that this translates into a $\wtilde O(d^3/\veps)$ computational complexity for finding an $\veps$-optimal point, answering an open question by Koren 2013. We first prove these new results for the simple case where the feasible set is a Euclidean ball. Then, to move to general convex sets, we use a reduction to Online Convex Optimization over the Euclidean ball. Our final algorithm for general convex sets can be viewed as a more computationally-efficient version of ONS.
	\end{abstract}
\arxiv{
	\clearpage 
	\tableofcontents
	\clearpage
}
	\section{Introduction}
	We consider the problem of Online Convex Optimization (OCO) with exp-concave loss functions. In this setting, an algorithm outputs vectors in a closed convex set $\K\subset \reals^d$ in rounds: At the beginning of round $t$, the algorithm outputs $\w_t\in \K$ based on the past history, then observes an exp-concave function $f_t\colon \K \rightarrow \reals_{\geq 0}$, which can be chosen by an adversary based on $\w_t$ and the history. The algorithm suffers loss $f_t(\w_t)$ and proceeds to the next round $t+1$. The performance of the algorithm is measured by its \emph{regret} against the best comparator vector $\w\in \K$ in hindsight after $T\in \mathbb{N}$ rounds: 
	\begin{align}
		\mathrm{Reg}_T \coloneqq 	 \sum_{t=1}^T f_t(\w_t) -  \min_{\w\in \K} \sum_{t=1}^T f_t(\w).\nn
	\end{align}
	Our goal in this paper is to design computationally-efficient algorithms that achieve (up to poly-log-factors) the optimal $O(d\ln T)$ regret for the online exp-concave setting \citep{mahdavi2015lower}.

	Many Machine Learning (ML) problems can be reduced to Online/Stochastic Exp-Concave Optimization. Of most relevance are certain \emph{online supervised learning} problems \citep{rakhlin2015online} that include, for example, Online Linear Regression (OLR), where at each round $t$, the algorithm receives a feature vector $\x_t\in \reals^d$ (that may be chosen adversarially), then outputs $\hat y_t\in \reals$ of the form $\hat y_t=\w_t^\top \x_t$ for some parameter vector $\w_t\in \K$ that is updated at every round. Then, an outcome $y_t$ is revealed and the algorithm suffers a loss $\ell(\hat y_t,y_t)$. The special case where $\ell$ is the square loss $\ell(\hat y_t,y_t)=(\hat y_t-y_t)^2$---where $\w\mapsto \ell(\w^\top\x,y)$ is exp-concave---has been extensively studied in the ML literature; see e.g.\,\citep{foster1991prediction,vovk1997competitive,vovk2001competitive,azoury2001relative,bartlett2015minimax,gaillard2019uniform}, and \cite[Chapter 11]{cesa2006prediction} for a thorough introduction to the topic. Other exp-concave losses that have also been extensively studied in the context of online learning include Cover's loss $\ell(\w,\x)=-\ln  (\w^\top \x)$ for the portfolio selection problem \citep{cover1991universal, luo2018efficient, mhammedi2022newton, jezequel2022efficient} and the logistic loss $\ell((\w,b),\x)= \ln (1+ \exp(-b \w^\top \x))$ for classification \citep{cover1991universal,foster2018logistic,agarwal2022efficient,MayoHE22}. Given the prevalence of ML problems that can be reduced to Online/Stochastic Exp-Concave Optimization, it is crucial to have efficient algorithms for the latter as echoed by the COLT open problem by \cite{koren2013open} that specifically asks for efficient algorithms in the stochastic setting.

	In constrained Online eXp-Concave Optimization (OXO) ensuring that the iterates $(\w_t)$ are within the feasible set $\K$, while guaranteeing the optimal $O(d\ln T)$ regret, typically requires a special type of projections, which constitute the main computational challenge behind the design of efficient OXO algorithms. For example, the Online Newton Step (ONS) \citep{hazan2007logarithmic}, one of the most popular algorithms for OXO, requires a \emph{Mahalanobis projection}, a.k.a.\,generalized projection, at each round $t$ where the iterate $\w_t$ steps outside of the set $\K$. Such projections typically require $\Omega(d^3)$ arithmetic operations (to perform an SVD decomposition) even for simple convex sets such as when $\K$ is a Euclidean ball. This means that in the stochastic setting, ONS-like algorithms can require $O(d^4/\veps)$ arithmetic operations to find an $\veps$-optimal point; see \citep{koren2013open} for details. \cite{koren2013open} asks for an algorithm that requires fewer than $O(d^4/\veps)$ arithmetic operations to find such a point.

	\paragraph{Contributions.}
	Our main contribution is a new, more efficient version of ONS, which does not require generalized projections thanks to the use of a self-concordant barrier. For the case where $\K$ is a Euclidean ball, our algorithm essentially outputs online Newton iterates with a log-barrier regularizer. The use of this barrier ensures that the iterates are always within $\K$ thanks to self-concordance properties of the barrier. The algorithm achieves the same regret bound as that of ONS (up to log-factors) and performs $\wtilde O(d^2)$ arithmetic operations per round, except for a $T^{-1/2}$ fraction of the rounds where the algorithm performs a matrix inversion. In contrast, standard ONS requires $\Omega(d^3)$ arithmetic operations (to perform an SVD decomposition) every round, in the worst-case. We are able to improve on the computational complexity of ONS by leveraging stability properties of the Newton iterates (which are conferred by our choice of regularizer) to efficiently approximate the inverse of certain Hessian matrices via Taylor expansions; for this reason, our approach falls under the category of \emph{Quasi-Newton Methods}\footnote{We note, however, that our new algorithm is distinct from the existing BFGS-based stochastic Quasi-Newton methods \citep{schraudolph2007stochastic, byrd2016stochastic,indrapriyadarsini2019stochastic,mokhtari2020stochastic}, which do not lead to an improved computational complexity for the stochastic exp-concave problem compared to ONS. } \citep{gill1972quasi}. 
	
	For general convex feasible sets, we use a reduction to OCO over the Euclidean ball. Using this reduction comes only with an additive $\wtilde O(T_{\sep}(\K))$ [resp.~$O(T_{\proj}(\K))$] computational cost per round for a centrally-symmetric [resp.~general] convex set $\K$, where $T_{\sep}(\K)$ denotes the cost of performing Separation [resp.~Euclidean projection] with respect to the set $\K$. In general, $T_{\proj}(\K)$ and $T_{\sep}(\K)$ can be much smaller than the cost of a Mahalanobis projection, which is required by ONS. We note that even projected Online Gradient Descent (OGD) requires $T_{\proj}(\K)$ arithmetic operations per round in the worst-case. We discuss the computational complexity of our algorithm in more detail in the sequel (Remark~\ref{rem:compcomplexity}).

	Finally, instantiating our OXO results in the stochastic exp-concave setting via a standard application of online-to-batch conversion (see e.g.\,\cite{cesa2006prediction}) leads to an algorithm that finds an $\veps$-optimal point using $\wtilde O(d^3/\veps)$ arithmetic operations. This answers one of the questions posed in the COLT open problem by \cite{koren2013open}. We conjecture that this number of arithmetic operations is the best one can hope for if one insists on a computational complexity that scales with $1/\veps$ (instead of $1/\veps^2$, for example). See App.~\ref{sec:regression} for more detail.

	\paragraph{Related Works.}
	The idea of using Newton steps with a barrier regularizer to build efficient online learning/optimization algorithms originated in \citep{abernethy2012interior}. There, the authors show that such online Newton iterates can approximate the Follow-The-Regularized-Leader (FTRL) iterates well enough to essentially inherit the regret guarantee of the latter. More recently, \cite{mhammedi2022newton} used online Newton iterates with a log-barrier for the simplex to approximate the iterates of the Mirror Descent-based algorithm BARRONS \citep{luo2018efficient} and build an efficient algorithm for the portfolio selection problem \citep{cover1991universal, van2020open}. However, the algorithms in \citep{abernethy2012interior, mhammedi2022newton} still require $\Omega(d^\omega)$ operations per round due to matrix inversion, where $\omega$ is the exponent of matrix multiplication. Reducing this computational cost to essentially $\wtilde O(d^2)$ per round is the main challenge we overcome in this paper by leveraging the stability of the Newton iterates when using the log-barrier for the Euclidean ball. 
	
	Our technique is inspired by one initially used in \citep{Vaidya1987} for efficiently solving linear programs by avoiding the computation of the full inverse of certain Hessian matrices at every iteration. However, their approach crucially relies on the feasible set being a polytope, and provides no computational advantage when the set is a Euclidean ball (we will reduce general OXO to OXO over a ball). Extending some of the ideas in \citep{Vaidya1987} to our setting is non-trivial and relies on recent results by \cite{mhammedi2022newton} on the stability of the Newton iterates with a particular choice of regularizer. We also note that our approach is different from previous ones that use, for example, sketching techniques to reduce the per-round computational complexity of ONS \citep{luo2016efficient}. Such techniques do not lead to a logarithmic regret in the OXO setting without additional assumptions.
	
	Extending our new efficient algorithm for OXO over the Euclidean ball to general convex sets by simply adapting the barrier to the set of interest fails because, beyond the Euclidean ball and polytopes, the barriers of other convex bodies are typically hard to compute. So, instead of taking this path, we leverage recent techniques in OCO \citep{cutkosky2018black,mhammedi2019,cutkosky2020parameter, mhammedi2022efficient} to reduce the exp-concave optimization problem to one over a Euclidean ball. In particular, we use the algorithm of \citep{mhammedi2022efficient} that reduces OCO over an arbitrary convex set to one over a ball for the purpose of efficient projection-free OCO\arxiv{\footnote{The algorithm of \cite{mhammedi2022efficient} is itself based on earlier algorithms by \cite{cutkosky2018black,cutkosky2020parameter} that were proposed for a slightly different purpose---that of reducing constrained OCO to unconstrained OCO.}}. Though we use the same algorithm as \cite{mhammedi2022efficient} for the reduction to OCO over a Euclidean ball, we need to extend their analysis to exp-concave losses; using their analysis directly leads to a $O(\sqrt{T})$ regret guarantee, which is sub-optimal in the OXO setting. We show that the reduction in \citep{mhammedi2022efficient} is naturally well suited to exp-concave losses, allowing our final algorithm to achieve near-optimal regret.

	\paragraph{Outline.}
	In Section~\ref{sec:prelim}, we present the notation and definitions we require in the main body of the paper. In Section~\ref{sec:newton}, we present our efficient algorithm for OXO over a Euclidean ball. In Section~\ref{sec:main}, we extend our results beyond the Euclidean ball using a reduction to online optimization on the latter. In Section \ref{sec:main}, we also instantiate our results in the stochastic exp-concave setting, where we answer one of the questions posed by \cite{koren2013open}.\arxiv{ We conclude with a discussion in Section~\ref{sec:discussion}.
	}

	\section{Preliminaries}
	\label{sec:prelim}
	Throughout, we let $\K$ be a closed convex subset of the Euclidean space $\reals^d$. We let $\|\cdot\|$ denote the Euclidean norm and $\cB(r)\subset \reals^d$ the Euclidean ball of radius $r>0$. We let $\gamma_{\cK}(\x)\coloneqq \inf \{ \lambda \geq 0\colon \x \in \lambda \cK\}$ be the \emph{Gauge function} of a convex set $\cK$, and $\cK^\circ\coloneqq \{\x \in \reals^d\colon  \inner{\x}{\y} \leq 1, \forall \y \in \cK \}$ be its \emph{polar} set  \citep{hiriart2004}. Further, we denote by $\inte \cK$ the interior of a set $\cK$. Our final algorithm requires access to either a Separation or a Euclidean projection Oracle for the set $\K$.
	\begin{definition}[Separation Oracle]
		A \emph{Separation Oracle} $\sep_{\cK}$ for a set $\cK$ is an Oracle that given $\u \in \reals^d$ either asserts that $\u\in \cK$ or returns $\w \in \reals^d$ such that $\inner{\w}{\u}> \inner{\w}{\v}$, for all $\v\in \cK$. We denote by $T_{\sep}(\cK)$ the computational complexity of one call to this Oracle. 
	\end{definition}

	\begin{definition}[Euclidean Projection Oracle]
		A \emph{Euclidean Oracle} $\proj_{\cK}$ for a set $\cK$ is an Oracle that given $\u \in \reals^d$ either asserts that $\u\in \cK$ or returns $\w \in \cK$ such that $\|\w -\u\|\leq\|\v-\u\|$, for all $\v\in \cK$. We denote by $T_{\proj}(\cK)$ the computational complexity of one call to this Oracle. 
	\end{definition}
	Our results can easily be extended to the case where only approximate Separation/Euclidean projection Oracles are available (see \citep{lee2018efficient} for definitions).

	This paper focuses on the online optimizing of exp-concave functions.
	\begin{definition}
		Let $\alpha>0$ and $\cK\subseteq \reals^d$ be a convex set. A function $f\colon \cK \rightarrow \reals$ is $\alpha$-exp-concave if $\x\mapsto e^{-\alpha f(\x)}$ is concave over $\cK$.
	\end{definition}
	For a twice differentiable function $f$, we denote by $\nabla^2 f(\u)$ [resp.~$\nabla^{-2} f(\u)$] its Hessian [resp.~inverse Hessian] at $\u$. We use the notation $\wtilde O(\cdot)$ to hide poly-log-factors in problem parameters. We denote by $\omega \in(2,3]$ the matrix multiplication constant. Our proofs use properties of self-concordant functions which we include in Appendix~\ref{sec:self-concordant}.
	
	\section{Efficient Online Exp-Concave Optimization Over a Ball}
	\label{sec:newton}
	In this section, we construct an efficient online algorithm (Alg.\,\ref{alg:pseudo}) for exp-concave optimization over the unit Euclidean ball $\cB(1)$. We later use Alg.\,\ref{alg:pseudo} as a subroutine in our main algorithm (Alg.\,\ref{alg:projectionfreewrapper}) for Online and Stochastic Exp-concave Optimization over general convex sets. The ``pseudocode'' (or the efficiently implementable version) of Algorithm \ref{alg:pseudo} is displayed in Algorithm \ref{alg:fastexpconcave} in the appendix. We carry out our analysis under the following Lipschitzness assumption on the sequence of losses.
	\begin{assumption} 
		\label{ass:grads}
		For $\Bfrak>0$, the sub-gradients $(\g_t)$ in Algorithm~\ref{alg:pseudo} satisfy $\|\g_t\|\leq \Bfrak,$ for all $t\geq 1$. 
	\end{assumption}
	The algorithm of this section, Alg.\,\ref{alg:pseudo}, essentially outputs approximate Newton iterates with respect to objective functions $(\Phi_t)$ that consist of the log-barrier for the unit Euclidean ball plus quadratic approximations of the observed losses $(\ell_t)$. In particular, for some parameters $(B, \eta, \beta)\in \reals^3_{>0}$, \arxiv{$(\Phi_t)$ are given by}
	\begin{gather}
		\Phi_t(\x) \coloneqq \Psi(\x) + \frac{d+B^2\eta}{2}\|\x\|^2 +\frac{\beta}{2}\sum_{s=1}^{t-1} \inner{\g_s}{\x-\w_{s}}^2 + \x^\top \sum_{s=1}^{t-1}\g_s, \label{eq:Phi} 
	\end{gather}
	where $\Psi(\x)\coloneqq  -\eta d\log (1-\|\x\|^2)$ and $(\g_t\in \partial \ell_t(\w_t))$ are the observed subgradients at the iterates $(\w_t)$ of Algorithm \ref{alg:pseudo}. With this, we show that the outputs $(\w_t)$ of Algorithm \ref{alg:pseudo} are approximate Newton iterates with respect to $(\Phi_t)$ in the sense that:
	\begin{align}
		\w_{t+1}\approx \w_t - \nabla^{-2}\Phi_{t+1}(\w_t)\nabla \Phi_{t+1}(\w_t),\quad \text{for all $t\in[T]$.} \label{eq:approx}
	\end{align}
	The regret analysis of Algorithm~\ref{alg:pseudo} then consists of showing that the Newton iterates with respect to $(\Phi_t)$ are good approximations of the FTRL iterates $(\x_t)$: \begin{align}\x_t\in \argmin_{\x\in \cB(1)}\Phi_t(\x)
		\label{eq:FTRL}
	\end{align} and that these FTRL iterates guarantee our target regret bound. What makes Algorithm~\ref{alg:pseudo} special is that it is able to efficiently approximate Newton iterates in the sense of \eqref{eq:approx}. To do this, the algorithm approximates the inverse of the Hessians $\nabla^2 \Phi_{t+1}(\w_t)$, for $t\geq 1$; for this reason, our approach falls under the category of Quasi-Newton Methods \citep{gill1972quasi}. Next, we discuss in more detail how Algorithm~\ref{alg:pseudo} is able to efficiently approximate Newton iterates in the sense of \eqref{eq:approx}.

		\begin{algorithm}
		\caption{OQNS: Online Quasi-Newton Steps Over the Euclidean ball. (Pseudocode in Alg.\,\ref{alg:fastexpconcave})}
		\label{alg:pseudo}
		\begin{algorithmic}[1]
			\REQUIRE Parameters $B,\eta,\beta, c >0$, and Taylor order $m$ as in Alg.\,\ref{alg:fastexpconcave}. \algcomment{$B,\beta$ needed for $(\Phi_t)$ in \eqref{eq:Phi}}
			\STATE Set $\u_1=\w_1 = \bm{0}$ and $A_0 = (2\eta d  + d +\eta B^2)^{-1} I$.
			\FOR{$t=1,\dots, T$}
			\STATE Play $\w_t$ and observe $\g_t \in \partial \ell_t(\w_t)$.
		\STATE Compute $A_{t} = \left(\frac{2 \eta  d I}{1-\|\u_t\|^2}  + (d + \eta B^2)I + \beta \sum_{s=1}^t \g_s \g_s^\top \right)^{-1}$.\label{line:expense}\algcomment{When $\u_t=\u_{t-1}$, $A_t^{-1}=A_{t-1}^{-1}+\beta \g_t\g_t^\top$, and so $A_{t}$ can be computed  in $O(d^2)$ using $A_{t-1}$ and a rank-one update of the inverse.} 
		\STATE Compute $H_t = \left(A_t^{-1}+\frac{4\eta  d \w_t \w_t^\top}{(1-\|\w_t\|^2)^2} \right)^{-1}$\algcomment{Computable in $O(d^2)$ with a rank-one update.}
			\STATE Define $\gamma_t= \frac{2\eta d}{1-\|\u_t\|^2}- \frac{2\eta d}{1-\|\w_t\|^2}$.\algcomment{The definition is such that $H^{-1}_t = \nabla^{2}\Phi_{t+1}(\w_t) +\gamma_t I$.}
			\STATE Define $\wtilde H_t = \sum_{k=1}^{m+1}\gamma_t^{k-1}  H_t^{k}$.  \algcomment{Approximates $\nabla^{-2}\Phi_{t+1}(\w_t)$ via $m^{\text{th}}$ order Taylor expansion}
			\STATE Compute $\w_{t+1}= \w_t - \wtilde H_t \nabla\Phi_{t+1}(\w_t)$. \label{line:taylor} \algcomment{Computable in $O(m d^2)$ (Lines \ref{line:set}-\ref{line:update} of Alg.\,\ref{alg:fastexpconcave})}
			\IF{$|\|\w_{t+1}\|^2-\|\u_{t}\|^2| \leq  c \cdot (1-\|\u_t\|^2)$}  \label{line:mark}
			\STATE Set $\u_{t+1}=\u_{t}$.\algcomment{No landmark update; $\u_t$ can be used for the next Taylor expansion}
			\ELSE 
			\STATE Set $\u_{t+1}=\w_{t+1}$. \algcomment{Updated landmark; this happens at most $\wtilde O(\sqrt{T})$ times (Lemma \ref{lem:movement})} \label{line:endmark}
		\ENDIF   
		\ENDFOR
	\end{algorithmic}
\end{algorithm}

	\subsection{Efficient Computation of Newton Iterates} Algorithm \ref{alg:pseudo} generates iterates $(\w_t)$ that satisfy the approximate equality in \eqref{eq:approx} without evaluating the inverse Hessian $\nabla^{-2}\Phi_{t+1}(\w_t)$ exactly at every round $t$. In particular, Algorithm \ref{alg:pseudo} computes the inverse Hessian $\nabla^{-2}\Phi_{t+1}(\w_t)$ exactly at most $\wtilde O(\sqrt{T})$ times after $T$ rounds, and approximates it using a Taylor expansion the rest of the time, resulting in a low amortized computational complexity. What makes this possible are certain self-concordance properties of the log-barrier $\Psi$ in the definition of $\Phi_t$. These properties imply that as long as the iterates $(\w_t)$ are stable enough, it is possible to efficiently approximate $\nabla^{-2}\Phi_{t+1}(\w_t)\nabla \Phi_{t+1}(\w_t)$ for most rounds $t\geq 1$ using Taylor expansions around a small number of landmark iterates $\w_{\tau_1},\dots, \w_{\tau_N}$, where $1\leq\tau_1\leq \dots \leq \tau_N \leq T$. At round $t\geq 1$, $\u_t$ in Algorithm \ref{alg:pseudo} represents the current landmark and $\u_{t}\neq \u_{t-1}$ only if $\w_{t}$ is ``far enough'' from the most-recent landmark $\u_{t-1}$, in which case $\w_{t}$ is set as the current landmark, i.e.\,$\u_{t}=\w_{t}$ (see Lines \ref{line:mark}-\ref{line:endmark} of Alg.\,\ref{alg:pseudo}). The next lemma, whose proof is in App.\,\ref{sec:movement_proof}, shows that the total number of unique landmarks used by Alg.\,\ref{alg:pseudo} is small relative to the number of rounds $T$. The proof of the lemma relies on the stability of the Newton iterates conferred by the non-linear terms in $(\Phi_t)$.
	\begin{lemma}[Stability]
		\label{lem:movement}
		Let $\beta,c \in(0,1)$, $B>0$, and $\eta \ge 1$. Further, let $(\u_t)$ be as in Algorithm \ref{alg:pseudo} with parameters $(B, \eta, \beta, c)$ and suppose that Assumption \ref{ass:grads} holds with $\Bfrak\leq B$. Then, it holds that
		\begin{align}
			\sum_{t=1}^{T-1}  \mathbb{I}\{ \u_{t+1} \neq \u_{t}\} \leq 8 \sqrt{\frac{2 T \ln (d+ B^2 T/d)}{c^2\eta\beta }}.\nn 
		\end{align}
	\end{lemma}
	On the rounds where the landmarks are updated (i.e.\,on the rounds $t$ where $\u_{t}\neq \u_{t-1}$) Algorithm~\ref{alg:pseudo} computes the inverse Hessian $\nabla^{-2}\Phi_{t+1}(\u_{t})$ exactly, which can be done in $O(d^{\omega})$ \citep{Cormen2009}. Next, we will show that on any other round (i.e.\,a round where $\u_{t}=\u_{t-1}$), Algorithm~\ref{alg:pseudo} efficiently approximates $\nabla^{-2}\Phi_{t+1}(\w_t)$ and computes $\w_{t+1}$ in $\wtilde{O} (d^2)$, implying (thanks to Lemma \ref{lem:movement}) a total computational complexity of at most $\wtilde{O}(d^2 T+ d^{\omega}\sqrt{T})$ for Algorithm \ref{alg:pseudo}.
	
\paragraph{Efficient approximation of the Hessian.}	On round $t\geq 1$, Algorithm~\ref{alg:pseudo}
	computes the iterate $\w_{t+1}$ via a Taylor expansion of order $m=O(\ln T)$ using the current landmark $\u_t$ as follows:
	\begin{align}
		\w_{t+1}=\w_t - \sum_{k=1}^{m+1} \gamma_t^{k-1} H_t^k \nabla \Phi_{t+1}(\w_t),\quad \text{where}\quad \gamma_t \coloneqq\frac{2 d\eta}{1-\|\u_t\|^2} - \frac{2d\eta}{1-\|\w_t\|^2}, \label{eq:caracter}
	\end{align}
and $H_t^{-1}= \nabla^{2}\Phi_{t+1}(\w_t) +\gamma_t I$. Now, with $\Theta\colon \gamma\mapsto (\nabla^2\Phi_{t+1}(\w_t)+ \gamma I)^{-1}$ and the definitions of $\gamma_t$ and $H_t$, the sum $\sum_{k=1}^{m+1} \gamma_t^{k-1} H_t^k$ is simply the $m$th-order Taylor expansion of $\Theta(0)=\nabla^{-2}\Phi_{t+1}(\w_t)$ at the point $\gamma_t$. This expansion provides an accurate approximation of $\Theta(0)=\nabla^{-2}\Phi_{t+1}(\w_t)$, when $\left|\frac{1-\|\u_t\|^2}{1-\|\w_t\|^2}-1\right|<1$ (i.e.~when $\gamma_t$ is small enough), which is an invariant of Algorithm~\ref{alg:pseudo} when $c<1$---see Line \ref{line:mark} of Algorithm \ref{alg:pseudo}. We formalize this in the next lemma.
	\begin{lemma}
		\label{lem:hessian}
		Let $\beta,c \in(0,1)$, $\eta \ge 1$, and $B>0$. Further, let $\gamma_t$ be as in \eqref{eq:caracter} and $(H_t)$ be as in Algorithm~\ref{alg:pseudo} with parameters $(B, \eta, \beta, c)$. Then, for any $m\geq 1$, we have
		\begin{align}
			\left\| \nabla^{-2}\Phi_{t+1}(\w_t) - \sum_{k=1}^{m+1} \gamma_t^{k-1} H_t^{k} \right\|	\leq  \frac{ c^m}{2\eta d\cdot(1-c)}. \nn 
		\end{align} 	
	\end{lemma}
One implication of this result is that the order $m$ of the expansion need only be of order $O(\ln T)$ (the exact choice of $m$ is specified in Alg.\,\ref{alg:fastexpconcave}) to obtain a $\text{poly}(\frac{1}{T})$-accurate approximation of the Hessian (which is all we need for our target regret). Thus, given the matrix $H_t$, the vector $\w_{t+1}$ in \eqref{eq:caracter} (i.e.~the output of Algorithm \ref{alg:pseudo}) can be computed in $O(md^2)=\wtilde O(d^2)$ using $m$ matrix-vector multiplications.

It remains to consider the computational cost of $H_t$ itself. First, note that $H_t$ we can be written as $H_t = \big(A_t^{-1}+\frac{4\eta d \w_t \w_t^\top}{(1 - \|\w_t\|^2)^2}\big)^{-1}$, where $A_t= \big(\frac{2 \eta  d I}{1-\|\u_t\|^2}  + (d + \eta B^2)I + \beta \sum_{s=1}^t \g_s \g_s^\top \big)^{-1}$. On the rounds where $\u_t\neq \u_{t-1}$, Algorithm \ref{alg:pseudo} computes $H_t$ using a full matrix inverse, costing $O(d^\omega)$. On the other hand, on a round $t$ where $\u_{t}=\u_{t-1}$, we have $A_{t}= (A_{t-1}^{-1}+ \beta \g_t \g_t^\top)^{-1}$. Thus, $A_t$ can be updated in $O(d^2)$ using $A_{t-1}$ (which is maintain by Algorithm~\ref{alg:pseudo}) and a rank-one update of the inverse. Further, since $H_t=\big(A_t^{-1}+\frac{4\eta d \w_t \w_t^\top}{(1 - \|\w_t\|^2)^2}\big)^{-1}$, $H_t$ can also be computed in $O(d^2)$ using $A_t$ and another rank-one update of the inverse. Since the number of rounds where $\u_{t}\neq \u_{t-1}$ is bounded by $\wtilde {O}(\sqrt{T})$ (Lemma~\ref{lem:close}), the total cost of computing $(H_t)$ and $(\w_{t})$ is at most $\wtilde{O}(d^2 T +d^\omega \sqrt{T})$.

	Having established that $(\w_t)$ approximate Newton iterates well (thanks to Lemma \ref{lem:hessian} and Line \ref{line:taylor} of Algorithm \ref{alg:pseudo}), we are now in a good position to bound the regret of Algorithm \ref{alg:pseudo}.
	\subsection{Regret Guarantee}
	\label{sec:surrogate}
 First, we note that when Assumption \ref{ass:grads} holds with $\Bfrak\leq B$, for some $B>0$, and the losses $(\ell_t)$ are $\alpha$-exp-concave, then for $\beta\leq \frac{1}{8 B}\wedge \frac{\alpha}{2}$, we have  
 \begin{align}
 	\ell_t(\w_t) - \ell_t(\w) \leq \inner{\w_t -\w}{\g_t}- \beta \inner{\w_t -\w}{\g_t}^2/2, \quad \text{for all }\w\in \cB(1). \label{eq:exp-concave}
 \end{align}
This result follows from \cite[Lemma 3]{hazan2007logarithmic}. Thus, to bound the regret of Algorithm~\ref{alg:pseudo}, it suffices to bound the sum $\sum_{t=1}^T (\inner{\w_t -\w}{\g_t}-\beta \inner{\w_t -\w}{\g_t}^2/2)$. To get a better handle on this sum, we will use the FTRL iterates $(\x_t)$ in \eqref{eq:FTRL}; in particular, we will add and subtract terms of the form $\inner{\x_t}{\g_t}$ and use H\"older's inequality to obtain the following bound (see details in Appendix~\ref{sec:omitted}):
\begin{align}
	&\sum_{t=1}^T \left(\inner{\w_t -\w}{\g_t}-\frac{\beta}{2} \inner{\w_t -\w}{\g_t}^2\right), \nn \\  &\leq 	\sum_{t=1}^T \left(\inner{\x_t -\w}{\g_t}-\frac{\beta}{2} \inner{\x_t -\w}{\g_t}^2\right) + \left(1 + 2\beta B  \right) \sum_{t=1}^T  \|\w_t -\x_t\|_{\nabla^2\Phi_t(\w_t)}  \|\g_t\|_{\nabla^{-2}\Phi_t(\w_t)}, \label{eq:red} 
	\end{align}
for all $\w\in \cB(1)$. The first sum on the RHS of \eqref{eq:red} is the regret of FTRL with respect to the surrogate losses $(\w\mapsto  \inner{\w}{\g_t} + \beta \inner{\x_t - \w}{\g_t}/2)$; these can be thought of as quadratic approximations of the actual losses $(\ell_t)$. The remaining term in \eqref{eq:red} can be bounded by $\sum_{t\in[T]}  \frac{1+2\beta B}{\sqrt{\eta}} \|\w_t-\x_t\|_{\nabla^2 \Phi_t(\w_t)}$ using that the local norms $(\|\g_t\|_{\nabla^{-2}\Phi_t(\w_t)})$ of the gradients are bounded by $1/\sqrt{\eta}$ thanks to our choice of regularizer $\Phi_t$ (we bound these local norms in Lemma~\ref{lem:gradsum} in the appendix). Thus, in light of \eqref{eq:red} and \eqref{eq:exp-concave}, to bound the regret of Algorithm \ref{alg:pseudo} it suffices to:\begin{itemize}
		\item[{\bf I.}] Bound the sum $\sum_{t\in[T]}   \|\w_t-\x_t\|_{\nabla^2 \Phi_t(\w_t)}$.
	\item[{\bf II.}] Bound the regret of FTRL with respect to the surrogate losses $(\w\mapsto  \inner{\w}{\g_t} + \beta \inner{\x_t - \w}{\g_t}/2)$.
	\end{itemize}
\paragraph{Point I: Bounding the sum of deviations.} The sum in Point I measures the deviation of the outputs of Algorithm \ref{alg:pseudo} from the FTRL iterates in the norms induces by the Hessians of the potentials $(\Phi_t)$. 
In Lemma \ref{lem:close} of App.~\ref{sec:helper}, we bound the sum in Point I from above by $\frac{16\sqrt{d}}{\beta\sqrt{\eta}} \ln (d + \frac{B^2T}{d})$ using that: 
\begin{itemize}
	\item[(i)] $(\w_t)$ are approximate Newton iterates; this follows from \eqref{eq:caracter} and Lemma \ref{lem:hessian} in the prequel.
	\item[(ii)]  The Newton iterates are close to the FTRL iterates $(\x_t)$ for an appropriate choice of parameters $(B,\eta, \beta)$. Our proof of the latter fact (see proof of Lemma \ref{lem:close}) is similar to one by \cite{mhammedi2022newton} who used damped Newton iterates to approximate FTRL iterates\arxiv{\footnote{The ability to approximate FTRL iterates with damped Newton steps has been leveraged before to design efficient online optimization algorithms \cite{abernethy2012interior,mhammedi2022newton}. In this paper, we use the `vanilla' Newton iterates instead of the damped ones.}}. 
	\end{itemize}
We will also use facts (i) and (ii) in the proof of our main theorem in this section (Theorem \ref{thm:regretNewton}) to show that the iterates $(\w_t)$ of Algorithm \ref{alg:pseudo} are always within $\cB(1)$.
	
\paragraph{Point II: Bounding the FTRL surrogate regret.}	It remains to bound the regret of FTRL with respect to the surrogate losses, which we do next (see proof in Appendix~\ref{sec:FTRL_proof}):
	\begin{lemma}[Surrogate regret of FTRL]
		\label{lem:FTRL}
		Let $\beta \in(0,1)$, $B>0$, and $\eta \ge 1$. If Assumption \ref{ass:grads} holds with $\Bfrak\leq B$, then the iterates $(\x_t)$ in \eqref{eq:FTRL} satisfy, for all $\w\in \inte\cB(1)$,
		\begin{align}
			\sum_{t=1}^T \left(\inner{\x_t -\w}{\g_t} - \frac{\beta}{2}\inner{\x_t -\w}{\g_t}^2\right) & \leq \Psi(\w) +\frac{d+\eta B^2}{2}\|\w\|^2  + \left(\frac{2d}{\beta}+\frac{32\sqrt{d}}{3\eta^{2}}\right) \ln (d + B^2 T/d).\nn 
		\end{align}
	\end{lemma}
We can now bound the (surrogate) regret of Alg.~\ref{alg:pseudo} using the bound on the sum of deviations in Point~I, the surrogate regret bound of FTRL, and \eqref{eq:red} (see proof in Appendix~\ref{sec:regretNewton_proof}):
	\begin{theorem}[Surrogate regret of Alg.~\ref{alg:pseudo}]
		\label{thm:regretNewton}
		Let $\beta\in(0,1/8)$, $c \in(0,1)$, $B>0$, and $\eta \ge 11$. Further, let $(\w_t)$ be the iterates of Algorithm \ref{alg:pseudo} with parameters $(B,\eta, \beta, c)$ and suppose that Assumption \ref{ass:grads} holds with $\mathfrak{B}\leq B$. Then, we have $\w_t \in \inte \cB(1)$, for all $t\geq 1$, and for all $\w\in \inte\cB(1)$,
		\begin{align}
			\sum_{t=1}^T \left(\inner{\w_t -\w}{\g_t}- \frac{\beta}{2} \inner{\w_t -\w}{\g_t}^2\right)\leq \Psi(\w) +\frac{d+\eta B^2}{2}\|\w\|^2  + \frac{(3 d+Bd^{\frac{1}{2}}) \ln (d + B^2 T/d)}{\beta}. \label{eq:actual}
		\end{align}
		The total computational complexity of the instance of Algorithm \ref{alg:pseudo} under consideration is bounded by \[
		O\left(m d^2  T  + c^{-1}d^\omega \sqrt{\beta^{-1} T \ln (d+ B^2 T/d)}\right),\] where $m=O(\ln T)$ is as in Algorithm \ref{alg:fastexpconcave}.
	\end{theorem}
\paragraph{From surrogate regret to actual regret.}
For any $\beta  \leq \frac{1}{8 B} \wedge \frac{\alpha}{2}$, Eq.~\eqref{eq:exp-concave} and Theorem \ref{thm:regretNewton} together immediately imply a $O(d\ln (d+T))$ bound on the actual regret $\sum_{t=1}^T (\ell_t(\w_t)-\ell_t(\w))$, for all $\w\in \cB(1-\frac{1}{T})$; shrinking the unit ball by $(1-\frac{1}{T})$ ensures that $\Psi(\w)$ in \eqref{eq:actual} is at most $O(\ln T)$. This guarantee can easily be extended to the whole unit ball $\cB(1)$ for $B$-Lipschitz losses $(\ell_t)$ using that $\ell_t((1-\frac{1}{T})\w)\leq \ell_t(\w)+  \frac{1}{T}(\ell_t(0)-\ell_t(\w))\leq \ell_t(\w)+  \frac{B}{T}$, for all $\w\in \cB(1)$, by convexity of the losses $(\ell_t)$.
Finally, we note that we state a bound on the surrogate regret in Theorem \ref{thm:regretNewton} instead of the actual regret because it will be convenient in the sequel when we generalize beyond the Euclidean ball.
	
\paragraph{Computational complexity.}Finally, we note that from Theorem \ref{thm:regretNewton}, the average per-iteration computational cost of Algorithm~\ref{alg:pseudo} is bounded by $\wtilde{O}(d^2 + d^\omega/\sqrt{T})$. We will later leverage this fact to show that Algorithm \ref{alg:pseudo} can be used in the stochastic exp-concave setting to find an $\veps$-optimal point in $\wtilde O(d^3/\veps)$ time complexity (this remains true even if we take $\omega =3$), which improves over the previous best $\wtilde O(d^4/\veps)$ \citep{koren2013open}. In the next section, we generalize the result of Theorem \ref{thm:regretNewton} beyond the Euclidean ball for both the online and stochastic settings.

	\section{Online/Stochastic Exp-Concave Optimization Over General Sets}
	\label{sec:main}
	In this section, we extend the results of Section~\ref{sec:newton} to a general convex set $\K$. One approach to achieving this would be to simply swap the log-barrier for the Euclidean ball in the definition of $\Phi_t$ in \eqref{eq:Phi} with a self-concordant barrier for the set $\K$. This would come with two main challenges. First, the barrier of a general convex set $\K$ that is, say, not a ball or a polytope, is typically hard to compute. This means that our algorithm from the previous section, which requires the gradients and Hessians of the barrier of the set of interest, is not a candidate for efficient exp-concave optimization over general convex sets. The second challenge is in the ability to approximate the inverse Hessian of an arbitrary barrier via a Taylor expansion. The fact that we were able to do this for the log-barrier of the Euclidean ball has to do with the special structure of the Hessian in this case. In particular, the Hessian of $\Psi$ (the log-barrier of the Euclidean ball) at a point $\w$ depends only on the norm $\|\w\|$ of $\w$ and the outer product $\w\w^\top$. If we ignore the outer product part, which is easy to deal with when it comes to computing the inverse Hessian thanks to the Sherman–Morrison formula, we are left with only a dependence in the norm $\|\w\|$. Therefore, a Taylor expansion in 1d is sufficient to approximate the inverse of the Hessian of $\Psi$ at $\w$. This is exactly what we do in Algorithm~\ref{alg:pseudo}. For the barrier of a general convex set, one would require a multivariate Taylor expansion to approximate the inverse of its Hessian, which can not always be done efficiently. 
	
	Given the challenges faced when changing the barrier regularizer to extend Algorithm \ref{alg:pseudo} to general convex sets, we instead reduce the OXO problem over general convex sets to one over the Euclidean ball. In the rest of this section, we present our new efficient OXO algorithm; Algorithm \ref{alg:projectionfreewrapper} with subroutine $\cA$ set as Algorithm \ref{alg:pseudo}. 
	Algorithm \ref{alg:projectionfreewrapper}, which is taken from \citep{cutkosky2020parameter,mhammedi2022efficient}, reduces OCO over any convex set $\K$ to OCO over a Euclidean ball. In fact, Algorithm~\ref{alg:projectionfreewrapper} (an algorithm over $\K$) essentially inherits the regret guarantee of its subroutine $\cA$, which is a subroutine over a Euclidean ball that contains $\K$. Thus, the problem becomes one of designing a subroutine $\cA$ with a good regret guarantee for exp-concave losses, which we have already tackled in Section \ref{sec:newton}. 
	Next, we describe and analyze Algorithm \ref{alg:projectionfreewrapper} in the general OXO setting before specializing the results to the stochastic exp-concave setting in Section~\ref{sec:stoch}. 
	\subsection{Efficient Online Exp-Concave Optimization via Reduction to the Ball}
	Before stating the regret guarantee of our algorithm in the online setting, we first formalize the assumptions we make starting with the Lipschitzness and exp-concavity of the losses.
	\begin{assumption} 
		\label{ass:exp-conc}
		For $\alpha>0$, the functions $(f_t\colon \K \rightarrow \reals)$ in Algorithm \ref{alg:projectionfreewrapper} are $\alpha$-exp-concave and $\sup\{\|\bm{\zeta}\| \colon \bzeta  \in \partial f_t(\w)\} \leq 1$, for any $\w \in \K$ \emph{(}i.e.\,$(f_t)$ are $1$-Lipschitz\emph{)}. Furthermore, $\K\subseteq \cB(1)$. 
	\end{assumption}
	We note that assuming that $(f_t)$ are $1$-Lipschitz and $\K\subseteq\cB(1)$ comes with no loss of generality as one can always re-scale the losses and the set $\K$ to satisfy this condition, provided the Lipschitz constant is known. When the Lipschitz constant is unknown, it is possible to adapt to it using known techniques such as those in \citep{mhammedi2019,cutkosky2019artificial,mhammedi2020}. For some of our results in this section, we will assume that the set $\K$ is centrally symmetric.
	\begin{assumption}
		\label{assum:symmetric}	
		The set $\K$ is centrally-symmetric, i.e.\,$\K= - \K$, and $\cB(1/\sqrt{d})\subseteq \K \subseteq \cB(1).$
	\end{assumption}
	Here again, there is no loss of generality in assuming that $\cB(1/\sqrt{d})\subseteq \K \subseteq \cB(1)$ when the set is centrally-symmetric since, in this case, it is always possible to apply a certain affine transformation (that puts the set into the isotropic position) to satisfy this condition (see e.g.\,\citep{lovasz2006simulated}). With our approach, we are able to leverage the fact that a set is centrally-symmetric for more efficient OXO (see Remark \ref{rem:compcomplexity} below), which is why we treat this case separately in what follows. 
	
	The next lemma essentially states that the instantaneous regret of Algorithm \ref{alg:projectionfreewrapper} can be bounded by that of  its subroutine $\cA$, and bounds the norm of the subgradients that the latter receives.
	\begin{algorithm}[t]
		\caption{A Reduction to Online Exp-Concave Optimization Over the Euclidean ball.}
		\label{alg:projectionfreewrapper}
		\begin{algorithmic}[1]
			\REQUIRE {\bf I)} An OCO algorithm $\cA$ over $\cB(1)\supseteq \K$; {\bf II)} A convex  function $\rho\colon \reals^d \rightarrow \reals_{\geq 0}\cup\{+\infty\}$. 
			\STATE Initialize $\cA$ and set $\w_{1} \in \cB(1)$ to $\cA$'s first output.
			\FOR{$t=1,2,\dots$}
			\STATE \label{line:S} Define $S(\w) = \inf_{\x \in \K} \rho(\w-\x)$. \algcomment{Easy to compute for the intended choices of $\rho$}
			\STATE \label{line:Sub} Set $\bnu_t \in \partial S(\w_t)$ and $\gamma_t  =\inner{\bnu_t}{\w_t}$.
			\STATE Play $\u_t=\mathbb{I}_{\{\gamma_t\geq 1\}} \w_t/\gamma_t +\mathbb{I}_{\{\gamma_t< 1\}}\w_t.$ \algcomment{$\u_t$ represents the "projection" of $\w_t$ onto $\K$.} 
			\STATE Observe subgradient $\bm{\zeta}_t \in \partial f_t(\u_t)$.
			\STATE  \label{line:2} Set $ \g_t = \bm{\zeta}_t -  \mathbb{I}_{\inner{\bm{\zeta}_t}{\w_t}<0}    \inner{\bm{\zeta}_t}{\u_t} \bnu_t$ 
			\STATE Set $\cA$'s $t$th loss function to $\ell_t : \w \mapsto \inner{{\g}_t}{\w}$.
			\STATE Set $\w_{t+1} \in \cB(1)$ to $\cA$'s $(t+1)$th output given the history $((\w_i,\ell_i)_{i\leq t})$.
			\ENDFOR
		\end{algorithmic}
	\end{algorithm}	
	
	\begin{lemma}
		\label{lem:inst}	
		Suppose that Assumption~\ref{ass:exp-conc} holds and let $\beta \leq \frac{1}{8}\wedge \frac{\alpha}{2}$. Then, the following holds:
		\begin{enumerate}
			\item[(a)] If $\rho(\cdot)$ in Algorithm~\ref{alg:projectionfreewrapper} is set to $\|\cdot\|$, then $\|\g_t\|\leq 1$.
			\item[(b)] If $\rho(\cdot)$ in Algorithm~\ref{alg:projectionfreewrapper} is set to $\gamma_{\K}(\cdot)$ and $\K$ satisfies Assumption~\ref{assum:symmetric}, then $\|\g_t\|\leq 1+\sqrt{d}$.
		\end{enumerate}
		Furthermore, the iterates $(\u_t)$ of the instance of Algorithm~\ref{alg:projectionfreewrapper} in either $(a)$ or $(b)$ satisfy
		\begin{align}
			\forall \w\in \K,  \ \ 	f_t(\u_t)-f_t(\w)\leq
			\inner{\bm{\zeta}_t}{\u_t- \w} - \frac{\beta}{2} 	\inner{\bm{\zeta}_t}{\u_t- \w}^2
			\leq 	\inner{\g_t}{\w_t- \w} - \frac{\beta}{2} 	\inner{\g_t}{\w_t- \w}^2.    \label{eq:monotone}
		\end{align}
	\end{lemma}
\arxiv{
	\begin{proof}[{\bf Proof}]
}
\colt{
		\begin{proof}
}		
[{\bf Case $(a)$}] When $\rho(\cdot)\equiv \|\cdot\|$, Alg.\,\ref{alg:projectionfreewrapper} matches \cite[Alg.\,1]{cutkosky2020parameter}, and so by \cite[Thm.\,2]{cutkosky2020parameter}, we have that I) $\|\g_t\|\leq \|\bzeta_t\|\leq 1$ (last inequality follows by Assump.\,\ref{ass:exp-conc}),  and II)
		\begin{align}
			\forall \w\in \K,\quad \inner{\bzeta_t}{\u_t-\w}\leq \inner{\g_t}{\w_t-\w}. \label{eq:monotone2}
		\end{align}
		This, together with the fact that the function $x\mapsto x - \beta x^2/2$ is non-decreasing over $[0,1/\beta]$ and $\inner{\g_t}{\w_t -\w}\leq \|\g_t\|\|\w_t-\w\|\leq 2\leq 1/\beta$ implies the second inequality in \eqref{eq:monotone}.
		
		[{\bf Case $(b)$}] Now, when $\rho(\cdot)\equiv \gamma_{\K}(\cdot)$ ($\gamma_{\K}(\cdot)$ is the gauge function of the set $\K$---see \S\ref{sec:prelim}), then Algorithm \ref{alg:projectionfreewrapper} matches \cite[Alg.\,1]{mhammedi2022efficient}, and so by \cite[Lemma~7]{mhammedi2022efficient}, we have $\|\g_t\|\leq (1+\kappa) \|\bzeta_t\|\leq 1+\kappa$, where $\kappa=R/r$ and $r,R$ are such that $\cB(r)\subseteq \K\subseteq \cB(R)$. By Assumption~\ref{assum:symmetric}, we have $\kappa=\sqrt{d}$ and so $\|\g_t\|\leq 1+\sqrt{d}$. On the other hand, since $\rho(\cdot)\equiv\gamma_{\K}(\cdot)$, the function $S$ in Algorithm~\ref{alg:projectionfreewrapper} satisfies $\inner{\bnu}{\w}=\gamma_{\K}(\w)$, for all $\bnu \in \partial S(\w)$ (see e.g.\,\cite[Lemma~6]{mhammedi2022efficient}). This means that $\gamma_t= \inner{\bnu_t}{\w_t}= \gamma_{\K}(\w_t)$, and so $\u_t = \mathbb{I}_{\{\gamma_{\K}(\w_t)\geq 1\}} \w_t/\gamma_{\K}(\w_t)  + \mathbb{I}_{\{\gamma_{\K}(\w_t)<1\}} \w_t$. Using this, and the triangle inequality, we get 
		\begin{align}
			|\inner{\g_t}{\w_t}|\leq  |\inner{\bzeta_t}{\w_t}| + |\inner{\bzeta_t}{\u_t} \inner{\bnu_t}{\w_t}|
			 =|\inner{\bzeta_t}{\w_t}| + |\inner{\bzeta_t}{\w_t} \inner{\bnu_t}{\u_t}| \stackrel{(*)}{\leq}2 |\inner{\bzeta_t}{\w_t}|\leq 2\leq \frac{1}{2\beta}, \nn 
		\end{align}
		where in $(*)$ we used that I) $\u_t\in \K$ and $\bnu_t\in \partial S(\w_t)\subseteq \K^\circ$ (see \cite[Lemma~6]{mhammedi2022efficient} for the set inclusion); and II) that $|\inner{\bnu}{\u}|\leq 1$, for all $\bnu\in \K^\circ$ and $\u \in \K$, which follows by definition of the polar set $\K^\circ$ (see \S\ref{sec:prelim}) and the fact that $\K$ is centrally-symmetric. By a similar argument, we also have $|\inner{\g_t}{\w}|\leq 1/(2\beta)$, for all $\w\in \K$, and so 
		\begin{align}
			|\inner{\g_t}{\w_t-\w}|\leq 1/\beta,\quad \forall \w\in \K.  \label{eq:bounde}
		\end{align}
	On the other hand, by \cite[Lemma 7]{mhammedi2022efficient}, we also have that $\inner{\bzeta_t}{\u_t-\w}\leq \inner{\g_t}{\w_t-\w}$, for all $\w\in \K$. Using this, Eq.\,\eqref{eq:bounde}, and that the function $x\mapsto x - \beta x^2/2$ is non-decreasing over $[0,1/\beta]$, implies the second inequality in \eqref{eq:monotone}. 
		
		Finally, for both cases $(a)$ and $(b)$, the first inequality in \eqref{eq:monotone} follows from Assumption \ref{ass:exp-conc} (i.e.\,the exp-concavity of the losses) and the range assumption on $\beta$ (see \cite[Lemma 3]{hazan2007logarithmic}).
	\end{proof}
	Using Lemma \ref{lem:inst}, we now state the main guarantee of Algorithm \ref{alg:projectionfreewrapper} when its subroutine $\cA$ is set to Algorithm \ref{alg:pseudo}, which outputs approximate Newton iterates over the unit Euclidean ball. Before stating this guarantee, we recall that $T_{\sep}(\K)$ [resp.~$T_{\proj}(\K)$] denotes the computational complexity of a Separation Oracle [resp.~Euclidean Projection Oracle] for the set $\K$. 
	\begin{theorem}
		\label{thm:main}
		Suppose that Assumption \ref{ass:exp-conc} holds and the subroutine $\cA$ in Algorithm \ref{alg:projectionfreewrapper} is set to Alg.\,\ref{alg:pseudo} with parameters $B>0$, $\eta \geq 11$, $c=1/4$, and $\beta \leq \frac{1}{8}\wedge \frac{\alpha}{2}$. Then, the following holds:
		\begin{itemize}
			\item[(a)] When $\rho(\cdot)$ in Alg.\,\ref{alg:projectionfreewrapper} is set to $\|\cdot\|$ and $B=1$, the regret of Alg.\,\ref{alg:projectionfreewrapper} after $T$ rounds is bounded by $O(d (\alpha^{-1} +G) \ln (d T))$, and the comp.~complexity is bounded by $\wtilde O((T_{\proj}(\K)+d^2)T + d^{\omega}T^{\frac{1}{2}})$.
			\item[(b)] When $\rho(\cdot)$ is set to $\gamma_{\K}(\cdot);$ $\K$ satisfies Assumption~\ref{assum:symmetric} (i.e.\,$\K$ is centrally-symmetric)$;$ and $B=1+\sqrt{d}$, the regret of Alg.\,\ref{alg:projectionfreewrapper} after $T$ rounds is bounded by $O(d (\alpha^{-1} +G) \ln (d T))$, and the total computational complexity is bounded by $\wtilde O((T_{\sep}(\K)+d^2)T + d^{\omega}T^{\frac{1}{2}})$.
		\end{itemize}
	\end{theorem}	
\arxiv{
	\begin{proof}[{\bf Proof}]
}
\colt{
	\begin{proof}
}	
We first analyze the regret then consider the computational complexity. Fix $\w \in \K$ and let $\wtilde \w \coloneqq (1-1/T)\w$. We bound the regret of Alg.\,\ref{alg:projectionfreewrapper} as 
		\begin{align}
			\sum_{t=1}^{T} (f_t(\u_t) - f_t(\w)) &= 		\sum_{t=1}^{T} (f_t(\u_t) - f_t(\wtilde \w))  +\sum_{t=1}^{T} (f_t(\wtilde \w) - f_t(\w)), \nn \\
			& \leq  \sum_{t=1}^{T} (f_t(\u_t) - f_t(\wtilde \w))  +1, \quad ((f_t) \text{ are }1\text{-Lipschitz and } \K\subseteq \cB(1))\nn \\
			& \leq  \sum_{t=1}^T\left( \inner{\w_t -\wtilde\w}{\g_t}- \frac{\beta}{2} \inner{\w_t -\wtilde \w}{\g_t}^2 \right)+1. \quad \text{(by Lemma \ref{lem:inst})}\label{eq:holds}
		\end{align} 
		The first sum on the RHS of \eqref{eq:holds} is the surrogate regret (see \S\ref{sec:surrogate}) of Alg.\,\ref{alg:pseudo}---the subroutine $\cA$ of Alg.\,\ref{alg:projectionfreewrapper}---against comparator $\wtilde \w$. To bound this surrogate regret, we will use Theorem~\ref{thm:regretNewton}. But first, we need to verify that Assumption~\ref{ass:grads}, under which Theorem \ref{thm:regretNewton} holds, is satisfied with $\Bfrak=B$ for the sequence $(\g_t)$. Thanks to Assumption~\ref{ass:exp-conc} [resp.~\ref{assum:symmetric}] and Lemma \ref{lem:inst}, Assumption~\ref{ass:grads} is satisfied for the sequence $(\g_t)$ with $\Bfrak=1$ [resp.~$\Bfrak=1+\sqrt{d}$] when $\rho(\cdot)\equiv\|\cdot\|$ [resp.~$\rho(\cdot)\equiv\gamma_{\K}(\cdot)$]. Thus, by Theorem~\ref{thm:regretNewton}, Eq.\,\eqref{eq:holds}, and the facts that $1- \|\wtilde \w\|^2\geq 1 - (1-\frac{1}{T})^2= \frac{2}{T}-\frac{1}{T^2}\geq \frac{1}{T}$ and $B \sqrt{d}\leq 2d$ (for both cases ($a$) and ($b$)), we get that in both cases of the theorem's statement:
		\begin{align}
			\sum_{t=1}^{T} (f_t(\u_t) - f_t(\w)) \leq \eta d \ln T +\frac{d+\eta d}{2}\|\wtilde\w\|^2  + \frac{5d \ln (d + T)}{\beta}. \label{eq:som}
		\end{align}
 Using that $\|\wtilde \w\|\leq 1$ in \eqref{eq:som} implies the desired regret bound.
		
The computational complexity of Algorithm~\ref{alg:projectionfreewrapper} is bounded by that of subroutine $\cA$, which by Lemma \ref{lem:movement} is less than $\wtilde O(d^2T+d^\omega\sqrt{T})$, plus $T$ times the computational complexity $T_{\grad}(S)$ of evaluating a subgradient of $S$ (this is required in Line \ref{line:Sub} of Alg.\,\ref{alg:projectionfreewrapper}). In case $(a)$, $S$ is differentiable everywhere except at the origin and $\nabla S(\w)=(\w- \Pi_{\K}(\w))/\|\w - \Pi_{\K}(\w)\|$, for $\w\in \reals^d \setminus \{\bm{0}\}$, where $\Pi_{\K}(\w)$ denotes the Euclidean projection of $\w$ onto $\K$. Thus, $T_{\grad}(S)\leq O(T_{\proj}(\K))$. In case $(b)$, i.e.\,when $\rho(\cdot)\equiv \gamma_{\K}(\cdot)$, we know that $T_{\grad}(S)$ can be bounded by $\wtilde {O}(T_{\sep}(\K))$ since $S$ and its subgradients can be evaluated using a Separation Oracle and a binary search (see e.g.\,\cite[\S{}B.2]{mhammedi2022efficient} and \citep{lee2018efficient}). 
	\end{proof}
	\begin{remark}[Computational Complexity] 
		\label{rem:compcomplexity} 
		The function $S$ is convex for the choices of $\rho$ in Theorem~\ref{thm:main}, and its subgradients (which are required in Alg.\,\ref{alg:projectionfreewrapper}) can be computed using either a Euclidean projection Oracle in case $\rho(\cdot)\equiv \|\cdot\|$, or a Separation Oracle with a binary search routine in case $\rho(\cdot)\equiv \gamma_{\K}(\cdot)$ \emph{(}see e.g.\,\citep{mhammedi2022efficient}\emph{)}. For many sets of interest $\K$, a Separation Oracle can be implemented in complexity $T_{\sep}(\K)\leq \wtilde O(d^2)$. A particularly relevant case is when a Membership Oracle for $\K$ can be implemented in $O(d)$ time, then $T_{\sep}(\K)\leq \wtilde {O} (d^2)$ \citep{lee2018efficient}. In many cases, we also have $T_{\proj}(\K)\leq \wtilde O(d^2)$ and, crucially, $T_{\proj}(\K)$ can be much smaller than the cost of a Mahalanobis projection, which is require by ONS. Finally, since $T_{\sep}(\K)\leq T_{\proj}(\K)$ in general, our approach is able to leverage that a set is centrally-symmetric for more efficient OXO.
	\end{remark}
	
	\subsection{Application to Stochastic Exp-Concave Optimization}
	\label{sec:stoch}
	We now instantiate our results in the stochastic setting where the sequence of losses $(f_t)$ are of the form $f_t(\cdot)=f(\cdot,\xi_t)$, where $(\xi_t)$ are i.i.d.~such that $f(\cdot)=\E[f(\cdot,\xi_t)]$ and $f$ is an exp-concave function. Specifically, we will make the following standard assumption in line with \cite{koren2013open}.
	\begin{assumption}
		\label{assump:stoch}
		The functions $(f_t)$ in Alg.\,\ref{alg:projectionfreewrapper} are such that $f_t(\cdot)=f(\cdot, \xi_t)$ and $\xi_1, \xi_2, \dots$ are i.i.d.~random variables in some set $\Xi$ and for all $\xi\in \Xi$,  $\w\rightarrow f(\w,\xi)$ is $\alpha$-exp-concave, for $\alpha>0$, and $\sup\{\|\bzeta\|\colon \bzeta \in \partial^{(1,0)} f(\w,\xi)\}\leq 1$. Furthermore, $\K \subseteq \cB(1)$.
	\end{assumption}
	The assumption that $\sup\{\|\bzeta\|\colon \bzeta \in \partial^{(1,0)} f(\w,\xi)\}\leq 1$ comes with no loss of generality as we can always re-scale the losses. We note that Assumption \ref{assump:stoch} implies Assumption \ref{ass:exp-conc}, and so the results of Theorem \ref{thm:main} apply. Using online-to-batch-conversion, the results of Theorem~\ref{thm:main} translate into excess-risk bounds in the stochastic setting\colt{ (the proof of the next theorem, which uses Theorem \ref{thm:main} and online-to-batch-conversion, is somewhat standard and we postpone it to Appendix~\ref{sec:mainstoch})}.
	\begin{theorem}
		\label{thm:mainstoch}
		Let $\veps\leq1/d$. Suppose that Assumption \ref{assump:stoch} holds and the subroutine $\cA$ in Algorithm \ref{alg:projectionfreewrapper} is set to Alg.\,\ref{alg:pseudo} with parameters $B>0, \eta \geq 11$, $c=1/4$, and $\beta = \frac{1}{8}\wedge \frac{\alpha}{2}$. Further, suppose that either 
		\begin{enumerate}
			\item[(a)] $\rho(\cdot)$ in Alg~\ref{alg:projectionfreewrapper} is set to $\|\cdot\|$ and $B =1$, or
			\item[(b)] $\rho(\cdot)$ is set to $\gamma_{\K}(\cdot)$, $\K$ satisfies Assumption~\ref{assum:symmetric} (i.e.\,$\K$ is centrally-symmetric), and $B =1+\sqrt{d}$.
		\end{enumerate}
		Then, for $T=\frac{d\ln (d/\veps)}{\alpha \veps}$ and $\bar \u_T\coloneqq \frac{1}{T}\sum_{t=1}^T \u_t$, where $(\u_t)$  are the iterates of Alg.\,\ref{alg:projectionfreewrapper}, we have
		\begin{align} 
			\E\left[f(\bar \u_T) - \inf_{\u\in \K}f(\u)\right]\leq O(\veps), \quad \text{where} \quad f(\cdot)\coloneqq \E[f(\cdot, \xi)].\nn
		\end{align}
	The comp.~costs in cases (a) and (b) are, respectively, $\wtilde O(\frac{d}{\veps}(d^2+T_{\proj}(\K)))$ and $\wtilde O(\frac{d}{\veps}(d^2+T_{\sep}(\K)))$.
	\end{theorem}	
	In the regime where $\veps>\frac{1}{d}$ (i.e.\,the regime not covered by Theorem~\ref{thm:mainstoch}), one can simply use projected OGD to find  an $\veps$-optimal solution with total computational complexity at most $O(\frac{1}{\veps^2}(d+ T_{\proj}(\K))) \leq O(\frac{d}{\veps}(d+ T_{\proj}(\K)))$ (where the last inequality follows by the fact that $\veps >\frac{1}{d}$), which is better than the computational complexity in Theorem \ref{thm:mainstoch} for general convex sets. (Though, we note that for small enough $\veps$, the complexity of OGD becomes worse than that of our algorithm and ONS.)
	 
We also note that the instance of Algorithm \ref{alg:projectionfreewrapper} in Theorem \ref{thm:mainstoch} can be used as a black box within an existing meta-algorithm due to \cite[Algorithm~1]{mehta2017fast} to achieve an excess risk guarantee with high probability (instead of in expectation). The computational complexity of the meta algorithm will only be worse than that of the instance of Alg.\,\ref{alg:projectionfreewrapper} in Theorem \ref{thm:mainstoch} by log factors in $T$ and $d$.  
	
	\paragraph{Implications for the open problem by \cite{koren2013open}.} The observation made in the previous paragraph and the results of Theorem \ref{thm:mainstoch}  directly answer one of the questions posed by \cite{koren2013open}. There, \cite{koren2013open} asks about the existence of an algorithm for stochastic exp-concave optimization over the Euclidean ball that can find an $\veps$-optimal point using fewer than $\wtilde O(d^4/\veps)$ arithmetic operations. For the case of a Euclidean ball, we have $T_{\proj}(\K)\leq O(d)$, and so the instance of Algorithm \ref{alg:projectionfreewrapper} in Theorem~\ref{thm:mainstoch} [resp.~OGD] finds an $\veps$-optimal point in less than $\wtilde O(d^3/\veps)$ time when $\veps \leq 1/d$ [resp.~$\veps>1/d$]. This remains true for general [resp.~centrally-symmetric] convex sets as long as $T_{\proj}(\K)\leq \wtilde O(d^2)$ [resp.~$T_{\sep}(\K)\leq \wtilde O(d^2)$], and even if we take the matrix multiplication constant to be $\omega=3$. We conjecture that $O(d^3/\veps)$ is the best one can do in general (see Appendix ~\ref{sec:regression}).

\arxiv{	
	\begin{proof}[{\bf Proof of Theorem \ref{thm:mainstoch}}]
		The result follows by Thm.\,\ref{thm:main} and standard online-to-batch conversion. If we let $\text{Reg}_T(\cdot)$ be the regret of Alg.\,\ref{alg:projectionfreewrapper} in response to the i.i.d.~loss functions $(f_t)$ and $\u_*\in \argmin_{\u\in \K} f(\u)$, where $f(\cdot)\coloneqq \E[f_t(\cdot)]$, then the average iterate $\bar \u_T$ of Alg.\,\ref{alg:projectionfreewrapper} after $T$ rounds satisfies
		\begin{align}
			\E[f(\bar \u_T)] - \inf_{\u \in \K} f(\u)\stackrel{(*)}{\leq}  \frac{1}{T}\sum_{t=1}^T \E[f_t(\u_t) -  f_t(\u_*)] = \frac{\E[\text{Reg}_T(\u_*)]}{T},   \label{eq:regretstoch}
		\end{align}
		where $(*)$ follows by Jensen's inequality and the fact that $\u_t$ is independent of $f_t$. Plugging the regret bounds of Algorithm~\ref{alg:projectionfreewrapper} in the settings of Theorem \ref{thm:main} implies that 
		\begin{align}
			\E[f(\bar \u_T)] - \inf_{\u \in \K} f(\u)\leq O(d(\alpha^{-1} +G) \ln (dT)/T). \label{eq:suboptimality}
		\end{align}
		By choosing $T = \frac{d}{\alpha \veps} \ln \frac{d}{\veps}$, we get that $\E[f(\bar \u_T)] - \inf_{\u \in \K} f(\u)\leq O(\veps)$. The bounds on the computational complexity follow directly from those in Theorem \ref{thm:main} by plugging-in the choice $T=\frac{d}{\alpha \veps} \ln \frac{d}{\veps}$ and using the fact that $\veps\leq 1/d$. This conclusion remains true even if we take $\omega =3$.
	\end{proof}
}

	\arxiv{
	\section{Discussion}
	\label{sec:discussion}
	In this paper, we presented a new efficient algorithm for Online Exp-Concave Optimization, which can be viewed as a more efficient version of the standard ONS algorithm. First, we designed an efficient algorithm over the Euclidean ball, then extend it to general convex sets thanks to a reduction to OCO over the former. When instantiating our results in the stochastic exp-concave setting, we obtain algorithms that can find an $\veps$-optimal point in $\wtilde O(d^3/\veps)$ time---improving over the previous best $\wtilde O(d^4/\veps)$. While this answers one question posed by \cite{koren2013open}, it still leaves one open---that of the existence of an algorithm that attains the optimal rate with only linear-time computation per iteration (our algorithm essentially requires $\wtilde O(d^2)$ computation per iteration). It is conceivable that, under some additional assumptions on the data-generating distribution, existing sketching techniques for efficient second-order learning, such as those used in \cite{luo2016efficient}, could be used along with our new algorithms to tackle this outstanding question. Without additional assumptions, we conjecture that it is not possible to do better than $O(d^3/\veps)$ if one insists on a computational complexity that scales with $1/\veps$ (instead of $1/\veps^2$, for example). One observation that lead us to this conjecture is that even in the simple special case of Linear Regression with the square loss (which is exp-concave), it is not clear if one can find an $\veps$-optimal point using fewer than $O(d^3/\veps)$ arithmetic operations (see Appendix \ref{sec:regression} for more detail).
}

	\clearpage
	
\arxiv{	
	\section*{Acknowledgment}
	ZM acknowledges support from the ONR through awards N00014-20-1-2336 and N00014-20-1-2394. We thank Adam Block and Ayush Sekhari for their helpful comments on the presentation. 
}
	
	\bibliography{biblio}

	\clearpage
	\appendix
	
	\colt{
		\tableofcontents
		\clearpage
	}
	
	\section*{Appendices}
	\addcontentsline{toc}{section}{Appendices}

	\section{Omitted Pseudocode and Regret Decomposition Details}
	\label{sec:omitted}
	\subsection{Pseudocode of Algorithm \ref{alg:pseudo}}
		\begin{algorithm}
		\caption{Pseudocode for OQNS (Algorithm \ref{alg:pseudo}).}
		\label{alg:fastexpconcave}
		\begin{algorithmic}[1]
			\REQUIRE Parameters $B,\eta,\beta, c >0$, and $m= \left\lceil- \log_c \left(\frac{ 12 (4+32/\eta^2)^2(2\eta d +  B^2 \eta + (B+2\beta B^2)T)^2 T}{(1-c)}\right)\right\rceil$.
			\STATE Set $\u_1=\w_1 = \bm{0}$, $V_0 = 0$, $A_0 =  I/(2 \eta d+d+ \eta B^2)$, and $\mathbf{S}_0=\G_0=\bm{0}$.
			\FOR{$t=1,\dots, T$}
			\STATE Play $\w_t$ and observe $\g_t \in \partial \ell_t(\w_t)$.
			\STATE \label{line:dil} Set $\G_t = \G_{t-1}+ \g_t$, $\mathbf{S}_t = \mathbf{S}_{t-1} +\g_t \g_t^\top \w_t$, and $V_t = V_{t-1}+\g_t\g_t^\top$. 
			\STATE \label{line:purge}Set $\bnabla_t=  \frac{2\eta d\w_{t}}{1-\|\w_t\|^2} +(d + \eta B^2 )\w_t +\beta V_t\w_t- \beta \mathbf{S}_t+ \G_t$. \algcomment{$\bnabla_t = \nabla \Phi_{t+1}(\w_t)$}
			\STATE \label{line:class}Set $A_t =   A_{t-1}  -  \frac{\beta A_{t-1} \g_t \g_t^\top A_{t-1}}{2+\beta \g_t^\top A_{t-1} \g_t}$ and $H_t = A_t - \frac{4\eta d  A_{t} \w_t \w_t^\top A_{t}}{(1-\|\w_t\|^2)^2+4 d\eta \w_t^\top A_{t} \w_t}$. 
			\STATE Set $\bm{\Delta}_{t}=\wtilde{\bm{\Delta}}_{t}=H_t \bnabla_t$. \label{line:set}
			\FOR{$k=1,\dots,m$}
			\STATE Update $\wtilde{\bm{\Delta}}_{t}\leftarrow   \left(\frac{2\eta d}{1-\|\u_t\|^2} - \frac{2\eta d }{1-\|\w_t\|^2}\right)  H_t  \wtilde{\bm{\Delta}}_{t}$.
			\STATE Update $\bm{\Delta}_{t} \leftarrow \bm{\Delta}_t+ \wtilde{\bm{\Delta}}_{t}$.
			\ENDFOR
			\STATE \label{line:update} Set $\w_{t+1}=\w_t -  \bm{\Delta}_t$. \algcomment{$\w_{t+1} \approx \w_{t}- \nabla^{-2} \Phi_{t+1}(\w_t) \nabla\Phi_{t+1}(\w_t).$} 
			\IF{$|\|\w_{t+1}\|^2-\|\u_{t}\|^2| \leq  c \cdot (1-\|\u_t\|^2)$} \label{eq:lineeight}
			\STATE Set $\u_{t+1}=\u_{t}$.
			\ELSE 
			\STATE Set $\u_{t+1}=\w_{t+1}$.
			\STATE \label{line:fbi} Set $A_t = \left(\frac{2\eta d I}{1-\|\w_{t+1}\|^2}  +dI +{\eta} B^2 I+{\beta} V_t\right)^{-1}$. \algcomment{$A_t= \left(\nabla^2 \Phi_{t+1}(\w_{t+1})- \frac{4\eta d\w_{t+1} \w_{t+1}^\top}{(1-\|\w_{t+1}\|^2)^2} \right)^{-1}$}  \label{line:endmark2}
			\ENDIF 
			\ENDFOR
		\end{algorithmic}
	\end{algorithm}

\subsection{Details on the Surrogate Regret Decomposition}
In this subsection, we provide more details on the regret decomposition in \eqref{eq:red} under Assumption \ref{ass:grads} with $\Bfrak\leq B$. For all $\w\in \cB(1)$, we have
\begin{align}
	&\sum_{t=1}^T \left(\inner{\w_t -\w}{\g_t}-\frac{\beta}{2} \inner{\w_t -\w}{\g_t}^2\right) \nn \\ 
		&= \sum_{t=1}^T \left(\inner{\x_t -\w}{\g_t}-\frac{\beta}{2} (\inner{\x_t -\w}{\g_t}+\inner{\w_t -\x_t}{\g_t})^2\right) +\sum_{t=1}^T \inner{\w_t -\x_t}{\g_t}, \nn \\ 
			&= \sum_{t=1}^T \left(\inner{\x_t -\w}{\g_t}-\frac{\beta}{2} \inner{\x_t -\w}{\g_t}^2\right) +\sum_{t=1}^T(1- {\beta}\inner{\x_t -\w}{\g_t}) \inner{\w_t -\x_t}{\g_t}- \frac{\beta}{2}\sum_{t=1}^T \inner{\w_t-\x_t}{\g_t}^2, \nn \\ 
			&\leq  \sum_{t=1}^T \left(\inner{\x_t -\w}{\g_t}-\frac{\beta}{2} \inner{\x_t -\w}{\g_t}^2\right) +\sum_{t=1}^T(1- {\beta}\inner{\x_t -\w}{\g_t}) \inner{\w_t -\x_t}{\g_t}, \nn \\ 
	 &\leq 	\sum_{t=1}^T \left(\inner{\x_t -\w}{\g_t}-\frac{\beta}{2} \inner{\x_t -\w}{\g_t}^2\right)  + \left(1 + 2\beta B  \right) \sum_{t=1}^T  \|\w_t -\x_t\|_{\nabla^2\Phi_t(\w_t)}  \|\g_t\|_{\nabla^{-2}\Phi_t(\w_t)}, \label{eq:surrogateregret}
\end{align}
where the last inequality follows by H\"older's inequality (to bound $|\inner{\w_t -\x_t}{\g_t}|$), and that $|\inner{\x_t-\w}{\g_t}|\leq 2 B$ (this follows by the fact that Assumption \ref{ass:grads} holds with $\Bfrak\leq B$ and that $\x_t,\w\in \cB(1)$).
	
	\section{Background on Self-Concordant Functions}
	\label{sec:self-concordant}
	In this section, we define self-concordant functions and present some of their properties that we make heavy use of in the proofs of our results. We start by the definition of a self-concordant function. For the rest of this section, we let $\cK$ be a convex compact set with non-empty interior $\inte \cK$. For a twice [resp.~thrice] differentiable function, we let $\nabla^2 f(\u)$ [resp.~$\nabla^3 f(\u)$] be the Hessian [resp.~third derivative tensor] of $f$ at $\u$. 
	\begin{definition}
		A convex function $f\colon \inte \cK \rightarrow \reals$ is called \emph{self-concordant} with constant $M_f\geq 0$, if $f$ is $C^3$ and satisfies {\bf I)} $f(\x_k)\to +\infty$ for $\x_k \to \x\in \partial \cK${\em ;} and {\bf II)}
		\begin{align}
			\forall \x\in  \inte\cK, \forall \u \in \reals^d, \quad |\nabla^3f(\x)[\u,\u,\u]| \leq 2 M_f \|\u\|^3_{\nabla^2 f(\x)}.\nn
		\end{align}
	\end{definition} 
	Note that by definition, if $f$ is self-concordant with constant $M_f\geq 0$ it is also self-concordant with any constant $M\geq M_f$. 
	For a self-concordant function $f$ and $\x\in \dom f$, the quantity $\lambda(\x, f)\coloneqq \|\nabla f(\x)\|_{\nabla^{-2}f(\x)}$, known as the \emph{Newton decrement}, will be instrumental in our proofs. The following two lemmas contain properties of the Newton decrement and Hessians of self-concordant functions, which we will use repeatedly throughout (see e.g. \cite{nemirovski2008interior,nesterov2018lectures}).
	\begin{lemma}
		\label{lem:properties}
		Let $f\colon \inte \cK\rightarrow \reals$ be a self-concordant function with constant $M_f\geq 1$. Further, let $\x\in \inte \cK$ and $\x_f\in \argmin_{\x\in \cK} f(\x)$. Then, {\bf I)} whenever $\lambda(\x,f)<1/M_f$, we have 
		\begin{align}
			\|\x -\x_f\|_{\nabla^2 f(\x_f)} 	\vee 	\|\x -\x_f\|_{\nabla^2 f(\x)} \leq {\lambda(\x,f)}/({1-M_f \lambda (\x,f)});  \nn
		\end{align}
		and {\bf II)} for any $M\geq M_f$, the \emph{Newton step} $\x_+\coloneqq \x - \nabla^{-2}f(\x)\nabla f(\x)$ satisfies $\x_+\in \inte \cK$ and 
		$ \lambda(\x_+,f)\leq M \lambda(\x,f)^2/ ( 1 - M \lambda(\x,f))^2.$
	\end{lemma}
	\begin{lemma}
		\label{lem:hessians}
		Let $f\colon \inte \cK\rightarrow \reals$ be a self-concordant function with constant $M_f$ and $\x \in \inte \cK$. Then, for any $\y$ such that $r\coloneqq \|\y - \x\|_{\nabla^2 f(\x)} < 1/M_f$, we have 
		\begin{align}
			(1-M_f r)^{2} \nabla^2 f(\y) \preceq \nabla^2 f(\x) \preceq (1-M_f r)^{-2}  \nabla^2 f(\x).\nn 
		\end{align}
	\end{lemma}
	The following result from \cite[Theorem 5.1.5]{nesterov2018lectures} will be useful to show that the iterates of our algorithms are always within the feasible set.
	\begin{lemma}
		\label{lem:deakin}
		Let $f\colon \inte \cK\rightarrow \reals$ be a self-concordant function with constant $M_f\geq 1$ and $\x \in \inte \cK$. Then, $\cE_{\x} \coloneqq \{\w \in \reals^d\colon \|\w-\x\|_{\nabla^{2}f(\x)}<1/M_f \}\subseteq \inte \cK$. Furthermore, for all $\w\in \cE_{\x}$, we have $$\|\w-\x\|_{\nabla^2 f(\w)} \leq \frac{\|\w- \x\|_{\nabla^2 f(\x)}}{1-M_f \|\w- \x\|_{\nabla^2 f(\x)}}.$$
	\end{lemma}
	Finally, we will also make use of the following result due to \cite{mhammedi2022newton}:
	\begin{lemma}
		\label{lem:inter0}
		Let $f\colon \inte \cK \rightarrow \reals$ be a self-concordant function with constant $M_{f}>0$. Then, for any $\x, \y \in \inte \cK$ such that $r\coloneqq \|\x-\y\|_{\nabla^2 f(\x)}<1/M_{f}$, we have 
		\begin{align}
			\|\nabla f(\x) - \nabla f(\y)\|^2_{\nabla^{-2}f(\x)}  \leq \frac{1}{(1-M_{f} r )^{2}} \|\y- \x\|^2_{\nabla^{2}f(\x)}.\nn
		\end{align}
	\end{lemma}
	We now have all the tools we require for the analysis of our OXO algorithms.

	\section{Proofs of Section \ref{sec:newton}}
	In this section, we prove the statements in Section \ref{sec:newton}. For this, we need a set of helper lemmas that we state in the next subsection. The proofs of the helper lemmas are in Appendix~\ref{sec:helper_proof}. 
	\subsection{Helper Lemmas}
	\label{sec:helper}
	First, we establish that the functions $(\Phi_t)$ are self-concordant.
	\begin{lemma}
		\label{lem:selfconcord}
		The function $\Phi_t$ in \eqref{eq:Phi} is a self-concordant function with constant $M_{\Phi_t}\leq 1/\sqrt{d\eta}$.
	\end{lemma}
	The next lemma gives a bound on the local gradients norms, which will be useful throughout (the proof is similar to ones in \citep{luo2018efficient,mhammedi2022newton}):
	\begin{lemma}
		\label{lem:gradsum}
		Let $\beta \in(0,1)$, $B>0$, and $\eta \ge 1$. Further, let $(\g_t)$ be such that $\|\g_t\|\leq B$, for all $t\geq 1$. Then, for any sequence $(\y_t) \subset \inte \cB(1)$, the potential functions $(\Phi_t)$ in \eqref{eq:Phi} satisfy 
		\begin{align}
			\forall t\in[T], \ \ \|\g_t\|^{2}_{\nabla^{-2}\Phi_{t}(\y_t)} \leq 1/\eta \quad \text{and} \quad 
			\sum_{t=1}^T \|\g_t\|^{2}_{\nabla^{-2}\Phi_{t}(\y_t)} \leq \frac{ d \ln (d+T B^2/d)}{\beta}.\nn
		\end{align}
	\end{lemma}
	The main technical heavy lifting in the paper is done in the proof of the next lemma. Some of the steps in the proof of the lemma that involve bounding the distance between $\w_t$ and $\x_t$ are similar to those found in \cite[Proof of Lemma 4.1]{abernethy2012interior} and \cite[Proof of Lemma 8]{mhammedi2022newton}. 
	\begin{lemma}[Master Lemma]
		\label{lem:close}
		Let $\beta,c \in(0,1)$, $B>0$, and $\eta \ge 1$. Further, let $(\w_t)$ be the iterates of Algorithm \ref{alg:pseudo} with parameters $(B, \eta, \beta, c)$ and suppose that Assumption \ref{ass:grads} holds with $\Bfrak\leq B$. Then, we have \emph{I)} $(\w_t)\subset  \inte \cB(1)$; and \emph{II)}
		\begin{gather}
			\forall t \geq 1, \ \ 	\frac{\sqrt{\eta d}}{4}\left(\|\w_t-\x_t  \|_{\nabla^2 \Phi_t(\w_t)}-\frac{1}{T}\right) \leq \frac{\sqrt{\eta d}}{2} \left(\lambda(\w_t, \Phi_{t})-\frac{1}{T}\right)\leq\lambda(\w_{t-1}, \Phi_{t})^2 \leq \frac{4}{\eta}.\nn
		\end{gather}
		Further, we have  $\sum_{t=1}^T \|\w_t - \x_t\|_{\nabla^2 \Phi_{t}(\w_t)}  \leq 1+ {\frac{16 \sqrt{d} }{3\beta \sqrt{\eta}}  \ln  (d+\frac{B^2T}{d})}$ and 
		\begin{align}
			\sum_{t=1}^T \|\w_t - \w_{t-1}\|_{\nabla^2 \Psi(\w_t)}^2 + \sum_{t=1}^T \|\w_t - \w_{t-1}\|_{\nabla^2 \Psi(\w_{t-1})}^2 \leq \frac{8 d\ln  (d+B^2T/d)}{\beta}. \label{eq:sum}
		\end{align}
	\end{lemma}

	\subsection{Proof of Lemma \ref{lem:movement}}
	\label{sec:movement_proof}
	For the proof of Lemma \ref{lem:movement}, we need the following elementary result.
	\begin{lemma}
		\label{lem:damped}
		Let $\Psi(\x)\coloneqq - \eta d \ln (1-\|\x\|^2)$. For any $\u, \w\in \cB(1)$, we have
		\begin{align}
			\frac{1}{\eta d} \|\w  - \u \|^2_{\nabla^2\Psi(\w )} \geq \frac{(\|\w \|^2-\|\u\|^2)^2}{(1-\|\w \|^2)^2}.\nn 
		\end{align}
	\end{lemma}
	
\arxiv{
	\begin{proof}[{\bf Proof}]
}
\colt{\begin{proof}
}
		Fix $\u,\w\in \cB(1)$. We have
		\begin{align}
			\frac{1}{2\eta d}	\|\w -\u\|^2_{\nabla^2 \Psi(\w )} & = (\w-\u)^\top \left(\frac{I}{1- \|\w\|^2} + \frac{2 \w\w^\top}{(1-\|\w\|^2)^2} \right) (\w-\u) ,\nn \\&=\frac{\|\w \|^2+ \|\u\|^2 - 2 \w ^\top \u  - 2 \|\w \|^2 \w ^\top \u +\|\w \|^4  +2(\w ^\top \u)^2 - \|\w \|^2 \|\u\|^2}{(1-\|\w\|^2)^2},\nn \\
			& = \frac{2(\|\w\|^2-\w ^\top\u)(1-\w ^\top\u)}{(1-\|\w \|^2)^2} +\frac{\|\u\|^2- \|\w \|^2}{1-\|\w \|^2},\nn \\
			& = \frac{2(\|\w\|^2-\w ^\top\u)(1-\|\w \|^2)}{(1-\|\w \|^2)^2} + \frac{2(\|\w\|^2-\w ^\top\u)^2}{(1-\|\w \|^2)^2}  +\frac{\|\u\|^2- \|\w \|^2}{1-\|\w \|^2}.\nn 
		\end{align}
		Now using that $-\w \u=2^{-1}(\|\w -\u\|^2 - \|\w \|^2 - \|\u\|^2)$, we get that 
		\begin{align}
			\frac{1}{2\eta d}	\|\w -\u\|^2_{\nabla^2 \Psi(\w )} 	& = \frac{\|\w- \u\|^2+\|\w \|^2-\|\u\|^2 }{1-\|\w \|^2} + \frac{2(\|\w\|^2-\w ^\top\u)^2}{(1-\|\w \|^2)^2}  +\frac{\|\u\|^2- \|\w \|^2}{1-\|\w \|^2},\nn \\
			& = \frac{\|\w- \u\|^2}{1-\|\w \|^2} + \frac{(\|\w- \u\|^2+\|\w \|^2-\|\u\|^2 )^2}{2(1-\|\w \|^2)^2}.\label{eq:delta} 
		\end{align}
		Now consider the function $f\colon X \rightarrow \frac{X}{1-\|\w \|^2} +\frac{(X+\|\w \|^2 -\|\u\|^2)^2}{2(1-\|\w \|^2)^2}$. Note that $\text{sgn}(f'(X))=\text{sgn}(X-\|\u\|^2+1)$. Thus, since $\|\u\|^2\leq 1$, the function $f$ is non-decreasing over $\reals_{\geq 0}$, and so $f(\|\w-\u\|^2)\geq f(0)$. Using this with \eqref{eq:delta}, we get 
		\begin{align}
			\frac{1}{\eta d}	\|\w -\u\|^2_{\nabla^2 \Psi(\w )} 	 \geq  \frac{(\|\w \|^2-\|\u\|^2)^2}{(1-\|\w \|^2)^2}.\nn
		\end{align}
	\end{proof}
	
	\arxiv{
	\begin{proof}[{\bf Proof of Lemma \ref{lem:movement}}]
	}
\colt{
		\begin{proof}{\bf of Lemma \ref{lem:movement}}
}
		Let $i_1, \dots, i_n$ be the rounds $t$ where $\u_t \neq \u_{t-1}$, and note that by Line \ref{line:mark} of Algorithm \ref{alg:pseudo}, we have 
		\begin{align}
			|\|\u_{i_{k+1}}\|^2 - \|\u_{i_k}\|^2| > c \cdot (1-\|\u_{i_k}\|^2), \quad \forall k \in[n-1].	\label{eq:sitone}
		\end{align}
		Further, let 
		\begin{align}
			\alpha_t \coloneqq \frac{\|\w_{t+1}\|^2 -\|\w_{t}\|^2}{1-\|\w_{t+1}\|^2}, \quad \text{and} \quad \mu_t \coloneqq \frac{\|\w_{t}\|^2 -\|\w_{t+1}\|^2}{1-\|\w_{t}\|^2}.\nn
		\end{align}
		Fix $k\in[n-1]$. Suppose that $(\sum_{t=i_{k}}^{i_{k+1}-1} \alpha_t )\vee (\sum_{t=i_{k}}^{i_{k+1}-1} \mu_t)\leq 1/2$ and let $m_k \coloneqq i_{k+1}-i_k$. In this case, by \eqref{eq:sitone} we have that 
		\begin{align}
			\ln(1+c)	&\leq  \left(\ln \frac{1-\|\u_{i_{k}}\|^2}{1-\|\u_{i_{k+1}}\|^2} \right) \vee  \left(\ln \frac{1-\|\u_{i_{k+1}}\|^2}{1-\|\u_{i_{k}}\|^2}\right), \nn \\ &\leq  \left(\ln \prod_{t=i_k}^{i_{k+1}-1} (1+\alpha_t)\right) \vee  \left(\ln \prod_{t=i_k}^{i_{k+1}-1} (1+\mu_t)\right)  ,\nn \\ &= \left(\sum_{t=i_k}^{i_{k+1}-1} \ln (1+\alpha_t)\right) \vee \left(\sum_{t=i_k}^{i_{k+1}-1} \ln (1+\mu_t)\right),\nn \\
			& \leq \ln \left(1+ \frac{1}{m_k}\sum_{t=i_k}^{i_{k+1}-1} \alpha_t  \right)^{m_k}\vee \ln \left(1+ \frac{1}{m_k}\sum_{t=i_k}^{i_{k+1}-1} \mu_t  \right)^{m_k}, \nn \quad \text{(Jensen)} \\
			& \leq \ln \left(1 +2 \sum_{t=i_k}^{i_{k+1}-1} \alpha_t\right)\vee \ln \left(1 +2 \sum_{t=i_k}^{i_{k+1}-1} \mu_t\right), \nn
		\end{align} 
		where the last inequality follows by the facts that $(\sum_{t=i_{k}}^{i_{k+1}-1} \alpha_t )\vee (\sum_{t=i_{k}}^{i_{k+1}-1} \mu_t)\leq 1/2$ and $(1+x)^r \leq 1+\frac{r x}{1-(r-1)x}$, for all $x\in(-1,\frac{1}{r-1}]$ and $r\geq 1$. Now, using that $\ln(1+x)\leq x$ for $x\geq 0$ and $\ln (1+x)\geq x/2$, for $x\in(0,1)$, we get that 
		\begin{align}
			\frac{c}{2}\leq\ln(1+c)& \leq  \left(2 \sum_{t=i_k}^{i_{k+1}-1} \alpha_t\right)\vee \left(2\sum_{t=i_k}^{i_{k+1}-1} \mu_t\right),\label{eq:new} \\
			& \leq  2\left(\sqrt{m_k \sum_{t=i_k}^{i_{k+1}-1} \alpha^2_t}\right)\vee \left(\sqrt{m_k\sum_{t=i_k}^{i_{k+1}-1} \mu^2_t}\right),\nn  \quad (\text{Jensen}) \\
			&  \leq  2 \sqrt{m_k  \sum_{t=i_k}^{i_{k+1}-1} \alpha^2_t +  m_k  \sum_{t=i_k}^{i_{k+1}-1} \mu^2_t}. 
			\label{eq:checkpoint}
		\end{align}
		So far, we have assumed that $(\sum_{t=i_{k}}^{i_{k+1}-1} \alpha_t )\vee (\sum_{t=i_{k}}^{i_{k+1}-1} \mu_t)\leq 1/2$. If this does not hold, then we have $(\sum_{t=i_{k}}^{i_{k+1}-1} \alpha_t )\vee (\sum_{t=i_{k}}^{i_{k+1}-1} \mu_t)\geq 1/2$. This implies \eqref{eq:new} from which \eqref{eq:checkpoint} follows. Now, \eqref{eq:checkpoint} implies 
		\begin{align}
			\sum_{t=i_k}^{i_{k+1}-1} \alpha^2_t +\sum_{t=i_k}^{i_{k+1}-1} \mu^2_t\geq \frac{c^2}{16 m_k}.\nn
		\end{align}
		Thus, by summing over $k=1,\dots,n-1$, and using Lemma \ref{lem:damped} and Lemma \ref{lem:close} (in particular \eqref{eq:sum}), we get 
		\begin{align}
			\frac{8  \ln (d+ B^2 T/d)}{\eta\beta } \geq \sum_{t=1}^T ( \alpha^2_t + \mu_t^2 )\geq \sum_{k=1}^n  \frac{c^2}{16 m_k} =\sum_{k=1}^n  \frac{c^2}{16 (i_{k+1}-i_k)}\geq \frac{c^2 n^2}{16 T},\nn 
		\end{align}
		where the last inequality follows by the fact that $x\mapsto 1/x$ is convex and Jensen's inequality.
		By rearranging, we get that 
		\begin{align}
			n \leq 8 \sqrt{\frac{2  T \ln (d+ B^2 T/d)}{c^2 \eta\beta }}.\nn
		\end{align}
	\end{proof}
	
	\subsection{Proof of Lemma \ref{lem:hessian}}
	\label{sec:hessian_proof}
	\arxiv{
	\begin{proof}[{\bf Proof}]
	}
\colt{
		\begin{proof}
}
		Fix $m\geq 1$ and let $\alpha_t \coloneqq \frac{ \|\w_t\|^2-\|\u_t\|^2}{1-\|\w_t\|^2}$. We have 
		\begin{align}
			\alpha_t = \frac{1-\|\u_t\|^2}{1-\|\w_t\|^2}- 1= -\frac{1-\|\u_t\|^2}{2\eta d} \gamma_t,\nn
		\end{align}
		where we recall that $\gamma_t =\frac{2 d\eta}{1-\|\u_t\|^2} - \frac{2d\eta}{1-\|\w_t\|^2}$.
		Note that $H_t$ in Alg.\,\ref{alg:pseudo} satisfies 
		\begin{align}
			H_t^{-1} &=  \frac{2 \eta d I}{1-\|\u_t\|^2} + \frac{4\eta  d\w_t \w_t^\top}{(1-\|\w_t\|^2)^2} + (d + \eta B^2)I + \beta V_t,\label{eq:slide} \\
			&  = \nabla^2 \Phi_{t+1}(\w_t) - \frac{2\eta d I}{1-\|\w_t\|^2}+\frac{2\eta  d I}{1-\|\u_t\|^2},\nn \\
			& = \nabla^2 \Phi_{t+1}(\w_t) - \frac{2\eta d I}{1-\|\u_t\|^2} \left( \frac{1-\|\u_t\|^2}{1-\|\w_t\|^2}-1\right),\nn \\
			& = \nabla^2 \Phi_{t+1}(\w_t) -  \frac{2\eta d \alpha_tI}{1-\|\u_t\|^2}. \nn 
		\end{align}
		Therefore, if we let $U_t \coloneqq (1-\|\u_t\|^2)H^{-1}_t/(2\eta d)$, we have 
		\begin{align}
			\nabla^{-2}\Phi_{t+1}(\w_t) & =  \left( \frac{2\eta d \alpha_t I}{1-\|\u_t\|^2} + H^{-1}_t\right)^{-1},\nn\\
			& = \frac{1-\|\u_t\|^2}{2\eta d} \left( \alpha_t I +  \frac{1-\|\u_t\|^2}{2\eta d}H^{-1}_t\right)^{-1},\nn \\
			& = \frac{1-\|\u_t\|^2}{2\eta d} ( \alpha_t I + U_t)^{-1},\nn \\
			& =    \frac{1-\|\u_t\|^2}{2\eta d}  U_t^{-1}  (I+\alpha_t U_t^{-1})^{-1}. \label{eq:cry} 
		\end{align}
		Now, by \eqref{eq:slide}, we have $U_t \succeq I$ and so $\|U^{-1}_t\| \leq 1$. Using this and that $|\alpha_t| \leq c <1$ (this is an invariant of Algorithm \ref{alg:pseudo}---see Line \ref{line:mark} of Alg.\,\ref{alg:pseudo}), we have 
		\begin{align}
			(1+\alpha_t U^{-1}_t)^{-1} =	 \sum_{k=0}^{\infty}(-\alpha_t)^k U_t^{-k}, \quad \text{and} \quad \left\| (1+\alpha_t U_t)^{-1} - \sum_{k=0}^{m} (-\alpha_t)^k U_t^{-k} \right\| \leq \frac{c^m}{1-c}. \nn 
		\end{align}
		Therefore, by \eqref{eq:cry} and the fact that $\|U_t^{-1}\|\leq 1$ we have 
		\begin{align}
			\left\| \nabla^{-2} \Phi_{t+1}(\w_t) -  \frac{1-\|\u_t\|^2}{2d\eta}  \sum_{k=1}^{m+1} (-\alpha_t)^{k-1} U_t^{-k} \right\|\leq \frac{(1-\|\u_t\|^2)\cdot c^m}{2\eta d \cdot(1-c)}.\nn
		\end{align}
		Now, the fact that $\frac{1-\|\u_t\|^2}{2d\eta}  \sum_{k=1}^{m+1} (-\alpha_t)^{k-1} U_t^{-k} =\sum_{k=1}^{m+1} \gamma_t^{k-1} H_t^{k}$ completes the proof.
	\end{proof}
	
	\subsection{Proof of Lemma \ref{lem:FTRL}}
	\label{sec:FTRL_proof}
\arxiv{
	\begin{proof}[{\bf Proof}]
	}
\colt{
		\begin{proof}
}
		Fix $\w\in \cB(1)$. Let $\phi_t(\x) \coloneqq \x^\top \g_t  + \beta  \inner{\g_t}{\x- \w_t}^2/2$ and $\phi_0(\x) \coloneqq \Psi(\x) + (d+\eta B^2)\|\x\|^2/2$, and note that $\Phi_t(\x)=\sum_{s=0}^{t-1}\phi_s(\x)$ and $\x_t \in \argmin_{\x\in \cB(1)} \sum_{s=0}^{t-1}\phi_s(\x)$. By \cite[Lemma 3.1]{cesa2006prediction}, we have
		\begin{gather}
			\sum_{t=0}^T \phi_t(\x_{t+1}) \leq \sum_{t=0}^T \phi_t(\w), \nn 
			\shortintertext{which implies that} 
			\sum_{t=1}^T \inner{\x_{t+1}-\w}{\g_t} \leq \Psi(\w) +\frac{d+\eta B^2}{2}\|\w\|^2 + \frac{\beta}{2} \sum_{t=1}^T \inner{\w_t -\w}{\g_t}^2.\label{eq:dude}
		\end{gather}
		Now, if suffices to bound the sum $\sum_{t=1}^T\inner{\x_t-\x_{t+1}}{\g_t}$. By Taylor's theorem, the exists $\y_t$ in the segment $[\x_t, \x_{t+1}]$ such that
		\begin{align}
			\Phi_{t+1}(\x_t) - \Phi_{t+1}(\x_{t+1}) &\geq   \nabla \Phi_{t+1}(\x_{t+1})^\top(\x_t - \x_{t+1}) +\frac{1}{2} \|\x_t-\x_{t+1}\|^2_{\nabla^2 \Phi_{t+1}(\y_{t})},\nn  \\
			& \geq  \frac{1}{2}\|\x_t-\x_{t+1}\|^2_{\nabla^2 \Phi_{t+1}(\y_{t})}, \label{eq:suppose}
		\end{align}
		where the last inequality uses the fact that $\x_{t+1}\in \argmin_{\x\in \cB(1)} \Phi_{t+1}(\x)$ is in the interior of $\cB(1)$ by self-concordance of $\Phi_{t+1}$.
		On the other hand, using the convexity of $\Phi_{t+1}$ and the fact that $\nabla \Phi_{t+1}(\x_t)=\nabla \phi_{t+1}(\x_t)+ \nabla \Phi_{t}(\x_t)=\nabla \phi_{t+1}(\x_{t})$ (by optimality of $\x_t$), we get that
		\begin{align*}
			\Phi_{t+1}(\x_t) - \Phi_{t+1}(\x_{t+1}) & \leq   \inner{\x_t-\x_{t+1}}{\nabla \phi_{t+1}(\x_t)}, \nn \\ & =   \inner{\x_t-\x_{t+1}}{\g_t}(1 + \beta \inner{\g_t}{\x_t - \w_t}),
			\nn \\ & \leq  \norm{\x_t - \x_{t+1}}_{\nabla^2\Phi_{t+1}(\y_{t})} \cdot \norm{\g_t}_{\nabla^{-2}\Phi_{t+1}(\y_{t})}\cdot (1 + \beta \inner{\g_t}{\x_t - \w_t}).
		\end{align*}
		Combining this and \eqref{eq:suppose}, we get \begin{align}
			\norm{\x_t - \x_{t+1}}_{\nabla^2\Phi_{t+1}(\y_{t})} \leq 2\norm{ \g_t}_{\nabla^{-2}\Phi_{t+1}(\y_{t})} (1 + \beta \inner{\g_t}{\x_t - \w_t}).\nn 
		\end{align}
		Using this and H\"older's inequality leads to 
		\begin{align}
			\inner{\g_t}{ \x_t - \x_{t+1}} \leq \norm{ \g_t}_{\nabla^{-2}\Phi_{t+1}(\y_{t})} \norm{\x_t - \x_{t+1}}_{\nabla^{2}\Phi_{t+1}(\y_{t})} &\leq 2\norm{ \g_t}_{\nabla^{-2}\Phi_{t+1}(\y_t)}^2(1 + \beta \inner{\g_t}{\x_t - \w_t}).\nn  
		\end{align}
		Thus, by summing this inequality for $t=1,\dots, T$, we get that
		\begin{align}
			\sum_{t=1}^T \inner{\g_t}{\x_t-\x_{t+1}}&\leq 	2\sum_{t=1}^T\norm{ \g_t}_{\nabla^{-2}\Phi_{t+1}(\y_t)}^2 + 2\beta \sum_{t=1}^T\norm{ \g_t}_{\nabla^{-2}\Phi_{t+1}(\y_t)}^2 \inner{\g_t}{\x_t - \w_t},  \nn \\& 
			\leq	2\sum_{t=1}^T\norm{ \g_t}_{\nabla^{-2}\Phi_{t+1}(\y_t)}^2  +2{\beta}\sum_{t=1}^T \norm{ \g_t}_{\nabla^{-2}\Phi_{t+1}(\y_t)}^2 \norm{ \g_t}_{\nabla^{-2}\Phi_{t}(\w_t)}\|\x_t - \w_t\|_{\nabla^2 \Phi_{t}(\w_t)},  \nn \\
			&\leq \frac{2 d \ln (d + B^2 T/d)}{\beta}+\frac{2\beta}{\eta^{3/2}}\sum_{t=1}^T \|\x_t - \w_t\|_{\nabla^2 \Phi_t(\w_t)},  \ (\text{Lem.\,\ref{lem:gradsum} and }\nabla^2 \Phi_{t+1}\succeq \nabla^2 \Phi_t)
			\nn \\
			& \leq \left(\frac{2d}{\beta}+\frac{32d^{1/2}}{3\eta^{2}}\right)\ln (d + B^2 T/d),\nn 
		\end{align}
		where the last inequality follows from the bound on $\sum_{t=1}^T \|\x_t - \w_t\|_{\nabla^2 \Phi_t(\w_t)}$ from Lemma \ref{lem:close}.
		Combining this with \eqref{eq:dude}, we get the desired bound.
	\end{proof}

	\subsection{Proof of Theorem \ref{thm:regretNewton}}
	\label{sec:regretNewton_proof}
	\arxiv{
	\begin{proof}[{\bf Proof}]
	}
	\colt{
	\begin{proof}
	}
		First, the fact that $(\w_t)\subset \inte \cB(1)$ follows from Lemma \ref{lem:close}.	Now, by the surrogate regret decomposition in \eqref{eq:surrogateregret} and the fact that $\|\g_t\|_{\nabla^2 \Phi_t(\w_t)}\leq 1/\sqrt{\eta}$ (Lemma \ref{lem:gradsum}), we have, for all $\w\in\inte \cB(1)$, 
		\begin{align}
& 	\sum_{t=1}^T \left(\inner{\w_t -\w}{\g_t}-\frac{\beta}{2} \inner{\w_t -\w}{\g_t}^2\right) \nn \\ 
&\leq 	\sum_{t=1}^T \left(\inner{\x_t -\w}{\g_t}-\frac{\beta}{2} \inner{\x_t -\w}{\g_t}^2\right)  + \left(1 + 2\beta B  \right) \sum_{t=1}^T  \|\w_t -\x_t\|_{\nabla^2\Phi_t(\w_t)} /\sqrt{\eta}. \label{eq:dec}
		\end{align}
	The first sum on the RHS \eqref{eq:dec} represents the surrogate regret of FTRL, and the second sum measures the deviation of the iterates of Alg.~\ref{alg:fastexpconcave} from the FTRL iterates $(\x_t)$.
Plugging the bound on the surrogate regret of FTRL [resp.~$\sum_{t=1}^T \|\w_t -\p_t\|_{\nabla^2\Phi_t(\w_t)}$] from Lemma \ref{lem:FTRL} [resp.~Lemma \ref{lem:close}] in \eqref{eq:dec}, we get that, for all $\w\in \inte  \cB(1)$, 
		\begin{align}
			&	\sum_{t=1}^T (\inner{\w_t -\w}{\g_t}- \beta \inner{\w_t -\w}{\g_t}^2/2) \nn \\ &\leq  \Psi(\w) +\frac{d+\eta B^2}{2}\|\w\|^2  + \left(\frac{2d}{\beta}+\frac{32\sqrt{d}}{3\eta^{2}}+ \frac{16\sqrt{d}}{3\beta \eta}+\frac{32 B \sqrt{d}}{3\eta}\right){\ln (d + B^2 T/d)},\nn \\
			&\leq   \Psi(\w) +\frac{d+\eta B^2}{2}\|\w\|^2  + \frac{(3 d+Bd^{1/2}) \ln (d + B^2 T/d)}{\beta},\label{eq:nurse}
		\end{align}
		where \eqref{eq:nurse} follows by the fact that $\frac{32\sqrt{d}}{3\eta^{2}}+ \frac{16\sqrt{d}}{3\beta \eta}\leq \frac{d}{\beta}$ and $32/(3\eta)\leq 1$ (since $\beta \in(0,1/8)$, $\eta \geq 11$, and $d\geq1$).
		
		We now look at the computational complexity of Algorithm \ref{alg:pseudo}. The most computationally expansive step in Algorithm \ref{alg:pseudo} is in Line \ref{line:expense}, which involves a full matrix inverse when $\u_{t}\neq \u_{t-1}$ (see also Line \ref{line:fbi} of Algorithm \ref{alg:fastexpconcave}; the pseudo-code of Alg.~\ref{alg:pseudo}). However, by Lemma~\ref{lem:movement}, the inverse need only be computed at most $O(c^{-1}\sqrt{\beta^{-1}{ T \ln (d+ B^2 T/d)}})$ times after $T$ rounds. The next most computationally expansive step in Alg.\,\ref{alg:pseudo} is in Line~\ref{line:taylor}, which involves computing the output $\w_{t+1}$. The output $\w_{t+1}$ in Line~\ref{line:taylor} can be computed in $O(m d^2)$ using $m$ matrix-vector multiplications (see Lines \ref{line:set}-\ref{line:update} of Alg.\,\ref{alg:fastexpconcave}). Thus, the claim on the total computational complexity of the algorithm follows by the fact that $m\leq O\left(1+\log_c\frac{d+\eta+\beta+B + T}{1-c}\right)$ (the exact choice of $m$ can be found in Algorithm \ref{alg:fastexpconcave}).
	\end{proof}

	\colt{
		\subsection{Proof of Theorem \ref{thm:mainstoch}}
		\label{sec:mainstoch}
			\begin{proof}
		The result follows by Thm.\,\ref{thm:main} and standard online-to-batch conversion. If we let $\text{Reg}_T(\cdot)$ be the regret of Alg.\,\ref{alg:projectionfreewrapper} in response to the i.i.d.~loss functions $(f_t)$ and $\u_*\in \argmin_{\u\in \K} f(\u)$, where $f(\cdot)\coloneqq \E[f_t(\cdot)]$, then the average iterate $\bar \u_T$ of Alg.\,\ref{alg:projectionfreewrapper} after $T$ rounds satisfies
		\begin{align}
			\E[f(\bar \u_T)] - \inf_{\u \in \K} f(\u)\stackrel{(*)}{\leq}  \frac{1}{T}\sum_{t=1}^T \E[f_t(\u_t) -  f_t(\u_*)] = \frac{\E[\text{Reg}_T(\u_*)]}{T},  \nn
		\end{align}
		where $(*)$ follows by Jensen's inequality and the fact that $\u_t$ is independent of $f_t$. Plugging the regret bounds of Algorithm~\ref{alg:projectionfreewrapper} in the settings of Theorem \ref{thm:main} implies that 
		\begin{align}
			\E[f(\bar \u_T)] - \inf_{\u \in \K} f(\u)\leq O(d(\alpha^{-1} +G) \ln (dT)/T). \nn
		\end{align}
		By choosing $T = \frac{d}{\alpha \veps} \ln \frac{d}{\veps}$, we get that $\E[f(\bar \u_T)] - \inf_{\u \in \K} f(\u)\leq O(\veps)$. The bounds on the computational complexity follow directly from those in Theorem \ref{thm:main} by plugging-in the choice $T=\frac{d}{\alpha \veps} \ln \frac{d}{\veps}$ and using the fact that $\veps\leq 1/d$. This conclusion remains true even if we take $\omega =3$.
	\end{proof}
}
	
	\section{Proofs of Helper Lemmas}
	\label{sec:helper_proof}
	
	\subsection{Proof of Lemma \ref{lem:selfconcord}}
	\arxiv{
	\begin{proof}[{\bf Proof}]	
	}
	\colt{
	\begin{proof}
	}
		First, we note that $\x\mapsto \Psi(\x)/(\eta d)=  -\ln (1- \|\x\|^2)$ is self-concordant with constant $1$ (see e.g.\,\cite[Exampled 5.1.1]{nesterov2018lectures}). Thus, $\Psi$ is a self-concordant function with constant $1/\sqrt{\eta d}$; this follows by the fact that if a function $f$ is self-concordant with constant $M_f$, then $\alpha f$, for $\alpha >0$, it is self-concordant with constant $1/\sqrt{\alpha}$ (see e.g.\,\cite[Corollary 5.1.3]{nesterov2018lectures}). On the other hand, since $\Phi_t(\x)$ is equal to $\Psi(\x)$ plus a quadratic in $\x$, then $\Phi_t$ is self-concordant with the same constant as $\Psi$ (see e.g.\,\cite[Corollary 5.1.2]{nesterov2018lectures}).
	\end{proof}
	
	\subsection{Proof of Lemma \ref{lem:gradsum}}
	\arxiv{
	\begin{proof}[{\bf Proof}]
	}
\colt{
\begin{proof}
}
		First note that $\eta \g_t\g_t^\top \preceq \eta B^2 I \preceq \nabla^2 \Phi_t(\y_t)-I$. This together with the fact that $\eta \ge \beta$ implies that
		\begin{gather}
			\|\g_t\|^2_{\nabla^{-2}\Phi_t(\y_t)} \leq \|\g_t\|^2_{(I + \eta \g_t \g_t^\top )^{-1}}= \g_t^\top(I + \eta \g_t \g_t^\top )^{-1}\g_t\leq 1/\eta,\nn \\ 
			\text{and} \qquad 	\|\g_t\|^{2}_{\nabla^{-2}\Phi_t(\y_t)}=  \g_t^\top\left( \nabla^2 \Psi(\y_t) +d I + \eta B^2 I+\beta V_{t-1} \right)^{-1} \g_t  \leq  \beta^{-1} \g_t^\top Q_{t}^{-1} \g_t. \label{eq:didit}
		\end{gather}
		where $Q_t \coloneqq d I  +\sum_{s=1}^t \g_s\g_s^\top$. Thus, by \eqref{eq:didit} and \cite[Lemma 11]{hazan2007logarithmic}, we have 
		\begin{align*}
			\|\g_t\|^{2}_{\nabla^{-2}\Phi_t(\y_t)}\leq 	\frac{1}{\beta} \sum_{t=1}^T \g_t^\top Q_t^{-1} \g_t \leq \frac{1}{\beta}\ln \frac{|Q_T|}{|Q_0|} \leq \frac{d\ln (d+TB^2/d)}{\beta}, 
		\end{align*}
		where the second inequality uses the fact $ \ln|Q_0|= d\ln d$ and by Jensen's inequality $ \ln|Q_T|\leq  d\ln \frac{\operatorname{tr} Q_T}{ d} \leq  d\ln \left( d+{\sum_{t=1}^T \norm{\g_t}^2_2 }/{ d}\right)\leq  d\ln( d+ B^2 T/d)$. 
		This completes the proof.
	\end{proof}
	
	\subsection{Proof of Lemma \ref{lem:close}}
	For the proof of Lemma \ref{lem:close}, we need two additional lemmas that we now state and prove:
	\begin{lemma}
		\label{lem:minPhi}
		Let $\eta,t \ge 1$. If $\w_1,\dots, \w_{t-1}\in \cB(1)$, then the FTRL iterate $\x_t$ in \eqref{eq:FTRL} satisfies:
		\begin{align}
			\frac{2\eta d}{1-\|\x_t\|^2} \leq 2(2\eta d +  B^2 \eta + (B+2\beta B^2)(t-1)).\nn 
		\end{align} 
	\end{lemma}
\arxiv{
	\begin{proof}[{\bf Proof}]
	}
\colt{
		\begin{proof}
}
		Since $\Psi(\x)$ is a self-concordant barrier, we have $\x_t \in \inte \cB(1)$. Thus, by the first-order optimality condition involving $\x_t$, we  have 
		\begin{align}
			\frac{2\eta d\x_t}{1-\|\x_t\|^2} + (d+B^2 \eta)\x_t + \beta \sum_{s=1}^{t-1} \g_s \g_s^\top (\x_t-\w_s) +\G_{t-1}=\bm{0}. \nn 
		\end{align}
		This implies that 
		\begin{align}
			\frac{2 \eta d\|\x_t\|}{1-\|\x_t\|^2} \leq d +  B^2 \eta + (B+2\beta B^2)(t-1). \label{eq:wehad}
		\end{align}
		If $\|\x_t\|\leq 1/2,$ then we are done since in this case $1/(1-\|\x_t\|^2)\leq 4/3\leq 2$. Otherwise, \eqref{eq:wehad} directly implies that 
		\begin{align}
			\frac{2\eta d}{1-\|\x_t\|^2}\leq 2(2\eta d +  B^2 \eta + (B+2\beta B^2)(t-1)).\nn 
		\end{align}
	\end{proof}
	
	\begin{lemma}
		\label{lem:bio}
		Let $\eta,t \ge 1$, $C>0$, and $\x_t \in \argmin_{\x\in \cB(1)} \Phi_t(\x)$. If $\w_1,\dots, \w_{t-1}\in \cB(1)$, then for any $\u \in \inte\cB(1)$ such that $\|\u - \x_t\|^2_{\nabla^2\Psi(\u)}\leq  \eta C^2$, we have
		\begin{align}
			\|\nabla \Phi_{t}(\u)\| &\leq  (4+2C)(2\eta d +  B^2 \eta + (B+2\beta B^2)(t-1)),\nn \\
			\text{and} \qquad   \nabla^2 \Phi_{t}(\u)  & \preceq  7(1+C)^2(2\eta d  +  B^2 \eta + (B+2\beta B^2)(t-1))^2  I.  \nn 
		\end{align}
	\end{lemma}
\arxiv{
	\begin{proof}[{\bf Proof}]
	}
\colt{\begin{proof}
}
		Fix $\u \in \inte\cB(1)$ such that $\|\u - \x_t\|^2_{\nabla^2\Psi(\u)}\leq  \eta C^2$. By Lemma \ref{lem:damped}, we have 
		\begin{align}
			C^2 & \geq 	\left(\frac{ \|\u\|^2 - \|\x_t\|^2}{1-\|\u\|^2}\right)^2 =  \left(\frac{ 1- \|\x_t\|^2}{1-\|\u\|^2} -1\right)^2.\nn 
		\end{align}
		This implies that
		\begin{align}
			\frac{2\eta d}{1-\|\u\|^2} \leq \frac{2\eta d \cdot(1+C)}{1-\|\x_t\|^2}\leq 2(1+C)(2\eta d  +  B^2 \eta + (B+2\beta B^2)(t-1)),  \label{eq:RGB}
		\end{align}
		where the last inequality follows by Lemma \ref{lem:minPhi}. Therefore, by the triangle inequality, we have
		\begin{align}
			\|\nabla \Phi_t(\u)\|&\leq \left\|\frac{2\eta d \u}{1-\|\u\|^2} \right\|+ \left\|(d+B^2 \eta)\u + \beta \sum_{s=1}^{t-1} \g_s \g_s^\top (\u-\w_s) +\G_{t-1} \right\| ,\nn\\  &\leq \frac{2\eta d}{1-\|\u\|^2}+(d + B^2 \eta) + 2\beta B^2 + B (t-1),\nn \\
			& \leq (4+2C)(2\eta d +  B^2 \eta + (B+2\beta B^2)(t-1)). \nn 
		\end{align}
		On the other hand, we have 
		\begin{align}
			\nabla^2\Phi_t(\u) & = \frac{2 \eta  d I}{1-\|\u\|^2} + \frac{4\eta  d \u \u^\top}{(1-\|\u\|^2)^2} + (d + \eta B^2)I + \beta V_t, \nn \\
			& \preceq 7(1+C)^2(2\eta d  +  B^2 \eta + (B+2\beta B^2)(t-1))^2  I.\nn 
		\end{align}
		where  the last inequality follows from \eqref{eq:RGB}.
	\end{proof}
	
	\arxiv{
	\begin{proof}[{\bf Proof of Lemma \ref{lem:close}.}] 
	}
\colt{	\begin{proof}{\bf of Lemma \ref{lem:close}}
}
		Define 
		\begin{gather}
			\wtilde\w_{t+1} \coloneqq \w_t - \nabla^{-2}\Phi_{t+1}(\w_t) \nabla \Phi_{t+1}(\w_t), \quad \text{and} \quad 	 \wtilde \bnabla_t \coloneqq \sum_{k=1}^{m+1}\left(\frac{2\eta}{1-\|\u_t\|^2} - \frac{2\eta}{1-\|\w_t\|^2}\right)^{k-1} H_t^{k} \bnabla_t,\nn 
		\end{gather}
		and note that $\w_{t+1}=\w_t-\wtilde \bnabla_t$. By induction, we will show that for all $s \geq 1$, 
		\begin{align}
			\w_s\in \mathrm{int}\  \cB(1)\  \& \ 
			\frac{\sqrt{\eta d}}{4}(\|\w_s-\x_s  \|_{\nabla^2 \Phi_s(\w_s)}-\epsilon) \leq \frac{\sqrt{\eta d}}{2} (\lambda(\w_s, \Phi_{s})-\epsilon)\leq \lambda(\w_{s-1}, \Phi_{s})^2 \leq \frac{4}{\eta}, \label{eq:induction} 
		\end{align}
		where $\epsilon =1/T$ and $\w_0 =\mathbf{0}$ by convention. The base case follows trivially since $\nabla \Phi_1(\w_0)=\nabla \Phi_1(\w_1)=\bm{0}$ and $\w_1 =\x_1$. Suppose that \eqref{eq:induction} holds for $s=t$. We will show that it holds for $s=t+1$. By the expression of $\Phi_{t+1}$ in \eqref{eq:Phi}, we have $\nabla \Phi_{t+1}(\w_t)=\g_t + \nabla \Phi_t(\w_t)$, and so by the fact that $(a+b)^2 \leq 2 a^2 +2 b^2$, we get 
		\begin{align}
			\lambda(\w_t, \Phi_{t+1})^2 &= \|\nabla \Phi_{t+1}(\w_t)\|^2_{\nabla^{-2}\Phi_{t+1}(\w_t)}  ,\nn \\ &\leq 2   \|\nabla \Phi_t(\w_t)\|^2_{\nabla^{-2}\Phi_{t}(\w_t)}  + 2 \|\g_t\|^2_{  \nabla^{-2}\Phi_{t}(\w_t)}, \ \ (\nabla^{2}\Phi_{t+1}(\cdot) \succeq \nabla^2 \Phi_t(\cdot)) \nn \\
			& =2   \lambda(\w_t, \Phi_t)^2  + 2 \|\g_t\|^2_{  \nabla^{-2}\Phi_{t}(\w_t)}, \label{eq:curz}\\
			& \le 2 (8^2/\eta^3 + 16 \epsilon/\eta^{3/2} + \epsilon^2)    + 2/\eta,\label{eq:penul}  \\ & \leq  4/\eta,
			\label{eq:postcurz} 
		\end{align}
		where in \eqref{eq:penul} we used the induction hypothesis in \eqref{eq:induction} for $s=t$ and the bound on $\|\g_t\|^2_{  \nabla^{-2}\Phi_{t}(\w_t)}$ from Lemma \ref{lem:gradsum}; and \eqref{eq:postcurz} uses the range assumptions on $\eta$ and that $\epsilon=1/T$. Since $\wtilde \w_{t+1}$ is the standard Newton step, Lemma~\ref{lem:properties} and the fact that $\lambda(\w_t,\Phi_{t+1})\leq (1-1/\sqrt{2})\sqrt{\eta d}$ (which follows from \eqref{eq:postcurz} and the range assumption on $\eta$), we have \begin{align}\lambda(\wtilde\w_{t+1},\Phi_{t+1}) \leq \frac{2}{\sqrt{\eta d}}\lambda (\w_t,\Phi_{t+1})^2.\label{eq:interm}
		\end{align} 
		Furthermore, since $\x_{t+1}$ is the minimizer of $\Phi_{t+1}$ and $\lambda(\wtilde \w_{t+1},\Phi_{t+1})\leq \sqrt{\eta d}/2$, we have $\|\wtilde\w_{t+1}- \x_{t+1}\|_{\nabla^2 \Phi_{t+1}(\wtilde\w_{t+1})}\leq 2 \lambda (\wtilde \w_{t+1}, \Phi_{t+1})$ (by Lemma~\ref{lem:properties} again). Combining this with \eqref{eq:interm} and \eqref{eq:postcurz} implies that  $\|\wtilde \w_{t+1}- \x_{t+1}\|^2_{\nabla^2 \Phi_{t+1}(\wtilde \w_{t+1})}\leq \eta C^2$ with $C= 16/\eta^2$. Thus, Lemma \ref{lem:bio} implies that 
		\begin{align}
			\label{eq:theprec}
			\nabla^2 \Phi_{t+1}(\wtilde \w_{t+1})   \preceq  7(1+32/\eta^2)^2(2\eta d +  B^2 \eta + (B+2\beta B^2)t)^2  I.
		\end{align}
		On the other hand, since $\bnabla_t = \nabla \Phi_{t+1}(\w_t)$ we have,
		\begin{align}
			\|\wtilde \w_{t+1}- \w_{t+1}\|&= \left\| \nabla^{-2}\Phi_{t+1}(\w_t) \nabla \Phi_{t+1}(\w_t)- \sum_{k=1}^{m+1}\left(\frac{2\eta}{1-\|\u_t\|^2} - \frac{2\eta}{1-\|\w_t\|^2}\right)^{k-1} H_t^{k} \bnabla_t \right\|,\nn \\ &=  \left\|\left( \nabla^{-2}\Phi_{t+1}(\w_t)  - \sum_{k=1}^{m+1}\left(\frac{2\eta}{1-\|\u_t\|^2} - \frac{2\eta}{1-\|\w_t\|^2}\right)^{k-1} H_t^{k}\right)  \nabla \Phi_{t+1}(\w_t)\right\|,\nn \\& \leq  \left\| \nabla^{-2}\Phi_{t+1}(\w_t)  - \sum_{k=1}^{m+1}\left(\frac{2\eta}{1-\|\u_t\|^2} - \frac{2\eta}{1-\|\w_t\|^2}\right)^{k-1} H_t^{k} \right\|\cdot\|\nabla \Phi_{t+1}(\w_t)\| ,\nn \\
			& \leq \frac{ c^m\cdot (4+32/\eta^2)(2\eta d +  B^2 \eta + (B+2\beta B^2)t)}{2\eta  d\cdot(1-c)}, \label{eq:yay}
		\end{align}
		where the last inequality follows by Lemma~\ref{lem:bio} (which holds with $C=16/\eta^2$ by \eqref{eq:induction}) and Lemma \ref{lem:hessian}. Combining \eqref{eq:yay} with \eqref{eq:theprec} implies that 
		\begin{align}
			\|\wtilde \w_{t+1}- \w_{t+1}\|_{\nabla^2\Phi_{t+1}(\wtilde \w_{t+1})} \leq \epsilon' \coloneqq 	\frac{ 3c^m\cdot (4+32/\eta^2)^2(2\eta d +  B^2 \eta + (B+2\beta B^2)t)^2}{2\eta d \cdot(1-c)}. \label{eq:bound}
		\end{align}
		We now show that this implies that $\w_{t+1} \in \cB(1)$. First, since $\w_t \in \cB(1)$ and 
		\begin{align}\|\wtilde \w_{t+1} -\w_t\|_{\nabla^2 \Phi_{t+1}(\w_t)}=\lambda(\w_t, \Phi_{t+1})\stackrel{(*)}{\leq} 2/\sqrt{\eta}< \sqrt{\eta d}/2,
			\label{eq:tilde}
		\end{align} where $(*)$ follows by \eqref{eq:postcurz}, we have that $\wtilde \w_{t+1}\in \cB(1)$ by Lemma~\ref{lem:deakin}. Now, by our choice of $m$ in Algorithm \ref{alg:pseudo}, we have $\epsilon' < \sqrt{\eta d}/4$, and so \eqref{eq:bound} implies that $\|\wtilde \w_{t+1}- \w_{t+1}\|_{\nabla^2\Phi_{t+1}(\wtilde \w_{t+1})} < \sqrt{\eta d}/4$. Therefore, $\w_{t+1}\in \cB(1)$ by Lemma \ref{lem:deakin}, since $\wtilde \w_{t+1}\in \cB(1)$.
		
		We now bound the Newton decrement $\lambda(\w_{t+1},\Phi_{t+1})$. First, by Lemma \ref{lem:deakin} and the fact that $\|\wtilde \w_{t+1}- \w_{t+1}\|_{\nabla^2\Phi_{t+1}(\wtilde \w_{t+1})} < \sqrt{\eta d}/4$, we have 
		\begin{align}
			\|\wtilde \w_{t+1}- \w_{t+1}\|_{\nabla^2\Phi_{t+1}(\w_{t+1})} \leq 2 \|\wtilde \w_{t+1}- \w_{t+1}\|_{\nabla^2\Phi_{t+1}(\wtilde \w_{t+1})}\leq 2\epsilon' <\sqrt{\eta d}/2. \label{eq:hessiancomp}
		\end{align}
		Using this, we get 
		\begin{align}
			\lambda(\w_{t+1},\Phi_{t+1})&= \|\nabla \Phi_{t}(\w_{t+1})\|_{\nabla^{-2}\Phi_{t+1}(\w_{t+1})}\nn \\ &  \leq \|\nabla \Phi_{t}(\wtilde \w_{t+1})\|_{\nabla^{-2}\Phi_{t+1}(\w_{t+1})} + \|\nabla \Phi_{t+1}(\w_{t+1})- \nabla \Phi_{t+1}(\wtilde\w_{t+1})\|_{\nabla^{-2}\Phi_{t+1}(\w_{t+1})} ,\nn \\
			& \leq (1-\epsilon'/\sqrt{\eta d})^{-1} \|\nabla \Phi_{t}(\wtilde \w_{t+1})\|_{\nabla^{-2}\Phi_{t+1}(\wtilde \w_{t+1})} + 2\|\wtilde \w_{t+1}- \w_{t+1}\|_{\nabla^2\Phi_{t+1}(\w_{t+1})},\label{eq:used} \\
			& \leq \lambda(\wtilde \w_{t+1}, \Phi_{t+1}) + 2\epsilon' \cdot \lambda(\wtilde \w_{t+1}, \Phi_{t+1})/\sqrt{\eta d} + 4 \epsilon', \nn \\
			& \leq \lambda(\wtilde \w_{t+1}, \Phi_{t+1}) +  \epsilon,  \label{eq:proev}
		\end{align}
		where \eqref{eq:used} uses Lemmas \ref{lem:hessians} and \ref{lem:inter0}, and the last inequality follows by \eqref{eq:interm}, \eqref{eq:postcurz} and the fact that $8 \epsilon' \leq \epsilon=1/T$ by the choice of $m$ in Algorithm \ref{alg:pseudo}. Now, since $\x_{t+1}$ is the minimizer of $\Phi_{t+1}$ and $\lambda(\w_{t+1},\Phi_{t+1})\leq \sqrt{\eta d}/2$ (by \eqref{eq:proev}, \eqref{eq:interm}, and \eqref{eq:postcurz}), we have $\|\w_{t+1}- \x_{t+1}\|_{\nabla^2 \Phi_{t+1}(\w_{t})}\leq 2 \lambda (\w_{t+1}, \Phi_{t+1})$ (by Lemma~\ref{lem:properties}). Combining this with \eqref{eq:proev}, \eqref{eq:interm}, and \eqref{eq:postcurz}, implies \eqref{eq:induction} for $s=t+1$, which concludes the induction. 
		
		We now use \eqref{eq:induction} together with \eqref{eq:curz} to bound the sums \begin{align}S_1 \coloneqq \sum_{t=1}^T \|\w_t - \x_t\|_{\nabla^2 \Phi_{t}(\w_t)}, \  S_2 \coloneqq \sum_{t=1}^T \|\w_t - \w_{t-1}\|_{\nabla^2 \Psi(\w_t)}^2, \  \& \   S_3 \coloneqq \sum_{t=1}^T \|\w_t - \w_{t-1}\|_{\nabla^2 \Psi(\w_{t-1})}^2. \nn
		\end{align} 
		To this end, we will first bound the sum $\sum_{t=1}^T \lambda (\w_t,\Phi_t)^i$, for $i=1,2$. Using that $\lambda(\w_{t+1},\Phi_{t+1}) \leq 2 (\eta d)^{-1/2}\lambda(\w_t,\Phi_{t+1})^2 +\epsilon$ (by \eqref{eq:induction}) and \eqref{eq:curz}, we get 
		\begin{align}
			\lambda(\w_{t+1}, \Phi_{t+1}) \leq 4 (\eta d)^{-1/2}	\lambda(\w_{t}, \Phi_{t})^2 + 4 (\eta d)^{-1/2} \|\g_t\|^{2}_{\nabla^{-2}\Phi_{t}(\w_t)}+\epsilon. \label{eq:match}
		\end{align}
		Summing \eqref{eq:match}, for $t=1,\dots, T$, rearranging, and using that $\lambda(\w_{T+1}, \Phi_{T+1})\geq 0$, we get 
		\begin{align}
			\sum_{t=2}^T \left(\lambda(\w_t, \Phi_t) - \frac{4}{\sqrt{\eta d}}  \lambda(\w_t, \Phi_t)^2\right) \leq   \frac{4}{\sqrt{\eta d}} \lambda(\w_1, \Phi_1)^2 + \frac{4}{\sqrt{\eta d}} \sum_{t=1}^T \|\g_t\|^{2}_{\nabla^{-2}\Phi_{t}(\w_t)}+T\epsilon.\nn 
		\end{align}
		Using \eqref{eq:induction} and the range assumption on $\eta$, we have $0\leq \frac{4}{\sqrt{\eta d}}\lambda(\w_t,\Phi_t)\leq  32/\eta^{2}+4\epsilon/\sqrt{\eta} \leq 1/4$. Therefore, we have
		\begin{align}
			\frac{3}{4}\sum_{t=1}^T \lambda (\w_t,\Phi_t) &\leq \lambda(\w_1,\Phi_1) + \frac{4}{ \sqrt{\eta d}}\sum_{t=1}^T \|\g_t\|^{2}_{\nabla^{-2}\Phi_{t}(\w_t)},\nn \\ &  \leq  \frac{1}{16} +  \frac{4}{ \sqrt{\eta d}} \sum_{t=1}^T \|\g_t\|^{2}_{\nabla^{-2}\Phi_{t}(\w_t)}\leq  \frac{1}{16} +  \frac{4 d \ln (d+B^2T/d)}{\beta\sqrt{\eta d}}, \label{eq:thesum}
		\end{align}
		where the last inequality follows by Lemma \ref{lem:gradsum} and the range assumption on $\eta$. Now, using the fact that $\x_t$ is the minimizer of $\Phi_t$, we have $\|\w_t-\x_t\|_{\nabla^{2}\Phi_t(\w_t)}\leq 2 \lambda(\w_t,\Phi_t)$, which implies the desired bound the sum $S_1$. We now bound $S_2$ and $S_3$. By Lemma \ref{lem:hessians} and \eqref{eq:hessiancomp}, we have 
		\begin{align}
			\|\w_{t+1}-\w_t\|_{\nabla^2\Psi(\w_{t+1})} & \leq  	2\|\w_{t+1}-\w_t\|_{\nabla^2\Psi(\wtilde\w_{t+1})}  , \nn \\ &\leq   	2\|\w_{t+1}-\wtilde \w_{t+1}    \|_{\nabla^2\Psi(\wtilde\w_{t+1})} + 2 \|\wtilde \w_{t+1}-\w_t\|_{\nabla^2 \Psi(\wtilde{\w}_{t+1})}, \nn \\
			& \leq 2\epsilon' +  \|\wtilde \w_{t+1}-\w_t\|_{\nabla^2 \Phi_{t+1}(\wtilde\w_{t+1})} \  \  \text{(by \eqref{eq:hessiancomp} and $\nabla^2 \Phi_{t+1}(\wtilde\w_{t+1})\succeq \nabla^2 \Psi(\wtilde\w_{t+1})$)},\nn \\
			& \leq \epsilon +  2\|\wtilde \w_{t+1}-\w_t\|_{\nabla^2 \Phi_{t+1}(\w_{t})}, \quad \text{(by \eqref{eq:tilde} and Lemma \ref{lem:deakin})} \nn \\
			& = \epsilon+\sqrt{2}\lambda(\w_t, \Phi_{t+1}). \quad  \text{(by \eqref{eq:tilde})} \label{eq:thenew}
		\end{align}
		On the other hand, since $\epsilon+\sqrt{2}\lambda(\w_t, \Phi_{t+1})\leq \sqrt{\eta d}/8$ (by \eqref{eq:induction}), Lemma \ref{lem:deakin} and \eqref{eq:thenew} imply that 
		\begin{align}
			\|\w_{t+1}-\w_t\|_{\nabla^2\Psi(\w_{t})}  \leq  \sqrt{2}\|\w_{t+1}-\w_t\|_{\nabla^2\Psi(\w_{t+1})}.\nn 
		\end{align}
		Now, to get the desired results, it suffices to bound the sum $S_2 \coloneqq \sum_{t=1}^T \|\w_t - \w_{t-1}\|_{\nabla^2 \Psi(\w_t)}^2$. Using \eqref{eq:curz}, \eqref{eq:thenew}, and the fact that $\lambda(\w_t,\Phi_{t+1})\leq 1$ (by \eqref{eq:induction}), we have 
		\begin{align}
			\sum_{t=1}^T \|\w_{t+1} - \w_{t}\|_{\nabla^2 \Psi(\w_{t+1})}^2 &\leq T\epsilon^2 +\sum_{t=1}^T 2^{3/2}  \epsilon  \lambda(\w_t,\Phi_{t+1}) +\sum_{t=1}^T 2\lambda(\w_t,\Phi_{t+1})^2, \nn \\
			& \leq 5T\epsilon^2 +   \sum_{t=1}^T \lambda(\w_t,\Phi_{t})^2/2+ 4  \sum_{t=1}^T\|\g_t\|^2_{\nabla^{-2}\Phi_t(\w_t)},\nn \\
			& \leq 5T\epsilon^2 +   \sum_{t=1}^T \lambda(\w_t,\Phi_{t})^2+ \frac{4 d \ln (d+B^2T/d)}{\beta}, \quad \text{(by Lemma \ref{lem:gradsum})}\nn \\
			& \leq 5T\epsilon^2 + \sum_{t=1}^T \lambda(\w_t,\Phi_{t})/2+ \frac{4 d \ln (d+B^2T/d)}{\beta}, \quad \text{($ \lambda(\w_t,\Phi_{t})\leq 1$ by \eqref{eq:induction})},  \nn \\
			& \leq 5T \epsilon^2  + \frac{6 d \ln (d+B^2T/d)}{\beta},\nn 
		\end{align}
		where the last inequality follows by \eqref{eq:thesum}.
	\end{proof}

	\section{Special Case of Linear Regression}
	\label{sec:regression}
	Without additional assumptions on the data-generating distribution, we conjecture that it is not possible to find an $\veps$-optimal point in Stochastic Exp-Concave Optimization using fewer than $O(d^3/\veps)$ arithmetic operations if one insists on a computational complexity that scales with $1/\veps$ (instead of $1/\veps^2$, for example). One observation that lead us to this conjecture is that even in the simple special case of Linear Regression (LR) with the square loss, it is not clear if one can find an $\veps$-optimal point using fewer than $O(d^3/\veps)$ arithmetic operations.
	
	In the LR setting with the square loss, one can assume that the covariates $\x_1,\x_2,\dots \in \reals^d$ are i.i.d., and $y_t = \w_\star^\top \x_t +\veps_t$, $t\geq 1$, for some fixed $\w_\star\in \reals^d$ (to be learned/approximated) and some i.i.d.~noise variables $\veps_1, \veps_2,\dots$. In this case, a natural estimator for $\w_\star$ is the Empirical Risk Minimizer (ERM) $\what \w \in \argmin_{\w\in \reals^d} \sum_{t=1}^T (\w^\top \x_t- y_t)^2$, which admits the closed form expression \[\what \w = (X^\top X)^\dagger X^\top \y,\] where $X$ [resp.~$\y$] is the matrix [resp.~vector] whose $t$th row is $\x_t^\top$ [resp.~$y_t$], and $\dagger$ denotes the pseudo-inverse. To ensure that $\what \w$ is an $\veps$-optimal point, in the sense that $\E[(\what \w^\top \x - \w_\star^\top \x)^2]\leq \veps$, standard generalization arguments say that $T$ needs be at least $\Omega (d/\veps)$, in general (see e.g.~\cite[Corollary 4.13]{philippe2019}). For such a $T$, $X$ is a matrix in $\reals^{d/\veps \times d}$, and so  evaluating even $X^\top X$ in the expression of $\what \w$ would require $d^3/\veps$ arithmetic operations. A similar number of arithmetic operations would, in general, be needed to project $\what \w$ onto a feasible set in case of constrained Linear Regression.

\end{document}